\documentclass{article}

\usepackage{arxiv}
\usepackage{mathmacros}
\usepackage[utf8]{inputenc} 
\usepackage[T1]{fontenc}    
\usepackage{url}            
\usepackage{booktabs}       
\usepackage{nicefrac}       
\usepackage{microtype}      
\usepackage{lipsum}		
\usepackage{graphicx}

\usepackage{xcolor}
\definecolor{mplblue}{RGB}{31, 119, 180}
\definecolor{mplorange}{RGB}{255, 127, 14}
\definecolor{mplgreen}{RGB}{44, 160, 44}

\usepackage{hyperref}
\hypersetup{
    colorlinks=true,
    linkcolor=mplblue,
    citecolor=mplgreen,
    urlcolor=mplorange
}

\usepackage{todonotes}
\let\oldtodo\todo
\renewcommand{\todo}[1]{\oldtodo[size=\footnotesize,color=mplorange!40]{#1}}

\usepackage[
    date=year,
    isbn=false,
    style=numeric,
    ]{biblatex}
\DeclareSourcemap{
  \maps[datatype=bibtex]{
    \map[overwrite]{
      \step[fieldsource=doi, final]
      \step[fieldset=url, null]
      \step[fieldset=urldate, null]
      \step[fieldset=eprint, null]
    }  
  }
}
\AtEveryBibitem{
  \clearfield{note}%
  \clearlist{language}
}

\usepackage{tikz}
\usetikzlibrary{
    shapes,arrows,calc,matrix,positioning,decorations.text,decorations.markings
    }

\usepackage{doi}

\usepackage{amsthm, thmtools}
\declaretheorem[name=Theorem,numberwithin=section]{thm}

\declaretheorem[name=Remark,style=remark,sibling=thm,qed=\qedsymbol]{rem}
\declaretheorem[name=Example,style=remark,sibling=thm,qed=\qedsymbol]{exa}

\newcommand{\orcid}{\raisebox{-1pt}{\includegraphics[scale=0.08]{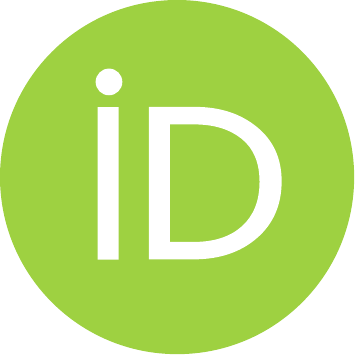}}\hspace{1mm}}

\title{Discovery of slow variables in a class of multiscale stochastic systems\\ via neural networks}
\renewcommand{\shorttitle}{Discovery of slow variables}

\author{
    \href{https://orcid.org/0000-0003-3317-8908}{\orcid}
    Przemys{\l}aw Zieli\'{n}ski
    \thanks{Corresponding author}
	\And
	\href{https://orcid.org/0000-0001-8074-1586}{\orcid}
	Jan S.~Hesthaven
	\And
	\textnormal{Institute of Mathematics, EPFL, Switzerland}\\[.5em]
	\texttt{\{przemyslaw.zielinski,jan.hesthaven\}@epfl.ch}
}

\hypersetup{
pdftitle={Discovery of slow variables},
pdfsubject={q-bio.NC, q-bio.QM},
pdfauthor={First Author, Second Author},
pdfkeywords={Multiscale dynamics, Dimensionality reduction, More},
}

\newcommand{\obs}{\mathcal{F}} 
\newcommand{\enc}{\mathcal{E}} 
\newcommand{\dec}{\mathcal{D}} 
\newcommand{\net}{\mathcal{N}} 
\newcommand{\proj}{\mathcal{P}} 
\newcommand{\smap}{\mathcal{S}} 

\newcommand{\levs}{\mathbb{L}} 
\newcommand{\ffib}{\mathbb{F}} 
\newcommand{\sman}{\mathbb{S}} 

\newcommand{\jac}[1]{J_{#1}}
\newcommand{\K}{\mathbb{K}}

\newcommand{\invm}{\mathsf{P}}
\newcommand{\invd}{\rh}

\begin{document}
\maketitle

\begin{abstract}
    Finding a reduction of complex, high-dimensional dynamics to its essential, low-dimensional ``heart'' remains a challenging yet necessary prerequisite for designing efficient numerical approaches. Machine learning methods have the potential to provide a general framework to automatically discover such representations. In this paper, we consider multiscale stochastic systems with local slow-fast time scale separation and propose a new method to encode in an artificial neural network a map that extracts the slow representation from the system. The architecture of the network consists of an encoder-decoder pair that we train in a supervised manner to learn the appropriate low-dimensional embedding in the bottleneck layer. We test the method on a number of examples that illustrate the ability to discover a correct slow representation. Moreover, we provide an~error measure to assess the quality of the embedding and demonstrate that pruning the network can pinpoint an essential coordinates of the system to build the slow representation.
\end{abstract}



\section{Introduction}
\label{sec:intro}
Extracting an effective dynamics from a high-dimensional systems remains one of the most challenging problems in computational modeling~\cite{froyland_computational_2014}. Nowadays, complex models with multiple time scales are found in a variety of domains, including bio-mechanics~\cite{bhattacharya_multiscale_2017, favino_multiscale_2018}, material research~\cite{praprotnik_multiscale_2008, van_der_giessen_roadmap_2020}, and climate science~\cite{majda_mathematical_2001, palmer_stochastic_2019}. The complexity of processes and scale separation in many such systems make the cost of direct simulation prohibitive. However, the long term evolution of certain aspects of this processes can often be described by much simpler \emph{reduced} dynamics that captures the essential behavior. Thus, knowing the correct aspects allows to build efficient model reduction techniques that not only decrease the dimension but also accelerate the simulations~\cite{debrabant_micro-macro_2017,katsoulakis_datadriven_2020,kevrekidis_equation-free_2009,legoll_effective_2010}.

To take on this challenge, we focus on multiscale stochastic systems with local slow-fast time scale separation and introduce a method to learn a map that extracts the slow representation from the system, which we term the \emph{slow variable}. We build the method on the approximation capabilities of neural networks, owing to their ability to deal, in principle, with very high-dimensional problems and the availability of easy computations of derivatives, which we employ to assess the accuracy of learned maps. Other approaches to this problem were explored based on eigenfunctions of the Koopman operator~\cite{froyland_computational_2014, froyland_trajectory-free_2016} or on manifold learning~\cite{singer_detecting_2009, dsilva_data-driven_2016}.  Besides systems with explicit slow-fast time separation, systems that exhibit \emph{metastable} behavior form another important class of multiscale models for which the problem of finding an appropriate reduced representation, in this context called \emph{reaction coordinate}, was addressed~\cite{mcgibbon_identification_2017,bittracher_transition_2018,bittracher_weak_2020}. Moreover, there exist data-driven techniques, like SINDY autoencoders~\cite{champion_data-driven_2019, champion_discovery_2019}, designed to discover the underlying governing equations in the reduced space.

To set the stage, consider the class of slow-fast stochastic systems which are Markovian models of two coupled equations: a slow and a fast SDE, where the time scale separation is given explicitly by a small parameter~$\ep$, see Eq.~\eqref{eq:sf_sde}. In such systems, the state variable $x$ decomposes into two sets of coordinates $(y,z)$, with $y$ containing the slow dynamics on the time scale of the order~$1/\ep$ and $z$ the fast fluctuations on the time scale of the order of~$\ep$. On the one hand, with the slow variable $y$ fixed, the fast process equilibrates rapidly to a conditional invariant distribution $\invm_{\!y}(dz)$ supported on the $y=\mrm{const}$ fiber. On the other hand, the evolution of the slow process can be effectively described on a certain submanifold of the state space $z=z(y)$ parametrized by the slow variable $y$. These two structures---the $y$-fibers and the $z(y)$ slow manifold---describe the slow-fast dynamics of the process and can be leveraged to design methods for effective simulation of slow-fast systems~\cite{givon_extracting_2004, e_analysis_2005}.

In this study, however, we target a more complicated situation than the prototypical slow-fast systems exhibit, in which the directions of the slow and fast dynamics do not align with the coordinates in which the system is defined~\cite{singer_detecting_2009, dsilva_data-driven_2016, froyland_trajectory-free_2016}. More precisely, we assume that the variables $(y,z)$ are observed through an unknown nonlinear function $f$ and the state variable reads $x=f(y,z)$. Therefore, the coordinates of the observed process in $x$ have both slow and fast dynamics mixed. Locally around a given $x$ there exist directions of both slow and fast variability, due to the noise in the system, and these directions depend on the current position in some nonlinear fashion, unlike in the slow-fast systems where we can globally split the coordinates into slow and fast. Moreover, since we assume the transformation $f$ is unknown, the slow dynamics of the process in $x$ remains hidden from the observer.

Our main objective is to encode in an artificial neural network a map, which we call the \emph{slow map}, from the observed coordinates $x$ to the slow variables $y$ without explicit knowledge of the nonlinear observation function~$f$. To achieve this, we develop a learning method that uses data from a number of simulations of the observed system. These data capture certain crucial yet computationally accessible features of the observed system. During training, the network explores the hidden noise correlations, contained in the dataset, to reveal the embedded lower-dimensional slow dynamics. This approach has potential to scale to high-dimensional systems and can naturally deal with vector valued slow maps.

To test the approach, we consider a class of stochastic differential equations which arise as a nonlinear transformation of a known slow-fast stochastic system. We term the underlying slow-fast SDE the \emph{hidden system} and the SDE resulting from the nonlinear transformation the \emph{observed system}. Since the observed system arises as a transformation of a slow-fast SDE, it shares certain properties with the hidden system. Crucial for our approach remains the fact that the state space of the observed system foliates into the fast fibers---a family of submanifolds on which the fast dynamics equilibrates rapidly---and contains an embedded slow manifold, transversal to the fibers, on which the long-term evolution happens. For the hidden slow-fast SDEs, the fast fibers always align with the $z$ coordinates, thus the slow and fast variables are clearly distinguishable. However, the geometry of these fibers in the observed state space becomes more complicated, rendering the distinction between slow and fast variables unclear.

The architecture of our networks consists of two stacked, fully connected feed-forward subnetworks: encoder and decoder, connected through the bottleneck layer of width equal to the dimension of the slow variables. Therefore, it is the same as for autoencoders; the difference lies in the way we use this architecture. The network is trained on a dataset consisting of representative points from a sample trajectory of the observed system together with their projections along the fast fibers onto the slow manifold, serving as the targets. In other words, we aim to approximate the projection map and we train the network in the supervised manner. 

Since the fully trained network represents the projection along the fast fibers, the encoder part learns a lower-dimensional representation of state space points that is constant along the fibers, i.e.~a slow map. We do not explicitly impose this feature on the encoder during the training. Rather, since the decoder learns a parametrization of the slow manifold that is connected to the projection along the fibers, the only continuous representation the encoder can converge to glues the points on the fibers.

We also address the issue of testing the encoder by introducing an error measure based on the orthogonality between the derivative of the encoder and the directions of local fast variability. A lack of orthogonality between the vectors spanning these two subspaces of the state space indicates that the level sets of the encoder deviate from the fast fibers of the observed system. Measuring the amount of this mismatch, we can assess the accuracy of the trained models. Since the derivatives of the encoder with respect to the input can be easily computed via backpropagation and the local fast directions via SVD, this yields a viable procedure to test the trained models.

Additionally, we demonstrate that applying a pruning technique during training, we can pinpoint the coordinates of the observed process that are essential for constructing the slow dynamics. It is often the case that in high-dimensional systems only a small number of system's coordinates is involved in the hidden slow dynamics. The remaining coordinates add to the ``nuisance'' dimensions related mostly to the noise in the system. Sifting out these inessential variables can help building efficient methods to simulate the long time dynamics of complex systems.

The paper is organized as follows. We begin in Section~\ref{sec:slow_observables} by examining the concept of a slow variable as a particular type of the observable of stochastic process. Next, we introduce in Section~\ref{sec:ms_class} a class of multiscale stochastic systems that we use to test our approach, which we describe in Section~\ref{sec:method}. Sections~\ref{sec:visu_2D} and~\ref{sec:error_prunning} contain numerical examples which illustrate and test the method.

The code for all numerical experiments in this paper was written in \href{https://www.python.org/}{Python}. For training the neural networks we use \href{https://pytorch.org/}{PyTorch} framework~\cite{paszke_pytorch_2019}. The code developed to train the networks and produce the figures can be found at \href{https://github.com/przemyslaw-zielinski/sf_nets}{github.com/przemyslaw-zielinski/sf\_nets}.

\section{Observables and slow variables of stochastic processes}
\label{sec:slow_observables}
In this section, we discuss in more detail the concept of a slow variable of multiscale processes and offer a rudimentary definition that encompasses the cases we study in this manuscript. Other attempts to elaborate on the definition of slow variables often use the spectral properties of the Koopman or transfer operator, see~\cite{froyland_computational_2014, froyland_trajectory-free_2016} or~\cite{bittracher_transition_2018} in the case of reaction coordinates for metastable systems. Here, we stick to stochastic terminology and anchor our description on the geometric decomposition of the state space and the relaxation properties of the fast dynamics. This reflects the intuition behind our method and the properties of slow-fast systems that are a point of departure for our approach.

Let $X_t$ be a stochastic process on $\R^D$. We call an \emph{observable} any continuous function $\obs\from\R^D\to\R^K$, where $K\leq D$, such that the expectation
\begin{equation*}
    t\mapsto\Exp[][\obs(X_t)]
\end{equation*}
exists for all $t\in[0,T]$, with $T\in(0,+\infty]$ a fixed final time of the simulation. Intuitively, an observable is a \emph{slow variable} if the image process on $\R^K$ through $\obs$, given by
\begin{equation*}
    Y_t = \obs(X_t),
\end{equation*}
evolves on a much slower time scale than $X_t$. To make this notion more precise, we introduce the following assumptions.

We suppose that $\obs$ is a $C^1$ function with full rank Jacobian matrix $\jac{\obs}$, and for each $y\in\R^K$ we define $\levs_y=\{x\in\R^D:\ \obs(x)=y\}$ as the~$y$-level set of~$\obs$. The family of all level sets $\{\levs_y:\ y\in\R^K\}$ forms a foliation of the state space $\R^D$. Consider the process $X_t$ conditioned to a fixed value $y$ of $\obs$
\begin{equation*}
    X_t\cbar\levs_y := \Exp[][X_t\cbar \obs(X_t)=y],\quad t\in[0,T].
\end{equation*}
We assume that for each $y$ there exists on $\levs_y$ an invariant probability measure $\invm_{\!y}$ of the process~$X_t\cbar\levs_y$.\medskip

\begin{rem}[On induced invariant measures]
If the process $X_t$ has itself an invariant measure $\invm=\invd\der{x}$ on~$\R^D$, the probability measures $\invm_{\!y}$ do exist~\cite{legoll_effective_2010}. Indeed, for each $y\in\R^K$ we can consider the probability measure $\invm_{\!y}$ such that
\begin{equation}
    \der{\invm_{\!y}} = \frac{1}{\Ga(y)}\rh(x)\det\big(\tp{\jac{\obs}(x)}\jac{\obs}(x)\big)^{-1/2}\der{\mathsf{S}_y},
\end{equation}
where $\mathsf{S}_y$ is the surface measure on $\levs_y$ (i.e., the Lebesgue measure on $\levs_y$ induced by the Lebesgue measure in the ambient Euclidean space $\R^D$) and $\Ga(y)$ is the normalization constant.
The co-area formula shows that if a random variable $X$ has law $\invm$, the law of $X$ conditioned to a fixed value $y$ of $\obs$ equals $\invm_{\!y}$.
\end{rem}

For any $x\in\R^D$, let $X^x_t = \Exp[][X_t\cbar X_0=x]$ and set $Y^x_t=\obs(X^x_t)$.
Intuitively, we call the process~$X_t$ multiscale and $\obs$ its slow observable, if starting at any $x\in\R^D$ there is a time window on which the process $X^x_t\cbar\levs_y$ already equilibrates to $\invm_{\!y}$ while the image process $Y^x_t$ evolves only slightly from its initial position $y$. Hence, we postulate an existence of a time scale $\tau$ such that the following relations hold approximately
\begin{gather}\label{eq:slowobs_condition}
\begin{aligned}
    X^x_\tau\cbar \levs_y &\sim \invm_{\!y},\\
    Y^x_\tau &\sim \de_y,
\end{aligned}
\end{gather}
where $\de_y$ is the Dirac mass at $y$. Relations~\eqref{eq:slowobs_condition} indicate that the slow observable $\obs$ collapses the fast fluctuations present in the dynamics of $X_t$ while retaining its slow evolution.

\begin{exa}[Slow and fast chemical reactions]

A chemical reaction network with $D$ species and $K$ reactions is described by the state vector $x\in\N_0^D$, containing the number of particles of each species, and instantaneous changes of the state due to transitions of the form
\begin{equation*}
    x\to x+\nu_k,
\end{equation*}
where $\nu_k$ is the stoichiometric vector describing the net change in the number of molecules of each species due to the $k$-th reaction~\cite{erban_stochastic_2019,winkelmann_stochastic_2020}. The rates for the reaction to take place are quantified by constants $\ga_k>0$ and the temporal evolution of the relevant process $X_t$ follows the dynamics in which the probability of the $k$-th reaction to occur is proportional to $\ga_k$. This dynamics can be simulated directly using the \emph{stochastic simulation algorithm} or by continuous approximation via \emph{chemical Fokker-Planck equation}. 

The reaction network is multiscale if we can decompose all rates $\ga_k$ into $D^f$ fast ones, indexed by $\K^f$, and $D^s$ slow ones, indexed by $\K^s$, such that $\ga_k\gg\ga_{k'}$ for all $k\in\K^f$ and $k'\in\K^s$. An observable $\obs\from\R^D\to\R^{D^s}$ of the~associated multiscale process $X_t$ is slow if it does not change during the fast reactions~\cite{e_nested_2007}, i.e., if for any~$x$
\begin{equation}\label{eq:slowobs_crn}
    \obs(x+\nu_k) = \obs(x),\quad\text{for all}\ k\in\K^f.
\end{equation}
To obtain a general representation of such observables it suffices to consider linear mappings satisfying~\eqref{eq:slowobs_crn}. For a linear map $\obs(x)=Bx$, where $B$ is a $D^s\times D$ matrix, to fulfill~\eqref{eq:slowobs_crn}, the rows $b_i$ of $B$ must satisfy $b_i\cdot\nu_k=0$, for all $k\in\K^f$. We can always find such a set $b_1,\dotsc,b_{D^s}$ of basis vectors in $\R^D$.

The level sets of linear slow variables are parallel hyperplanes of dimension $D^f$ that do not typically align with the coordinates of the state space, due to the mixing of species during reactions. Thus multiscale chemical reactions exemplify the generalization of the slow-fast systems we consider in this paper. However, since we focus on nonlinear mixing of slow and fast directions in the process $X_t$, the linear form of slow variables makes it too simple a case for our considerations. We refer to Section~\ref{sec:example} for a more pertinent example.
\end{exa}

\section{A class of multiscale systems}\label{sec:ms_class}
In this section, we introduce a particular class of models that we consider in this manuscript. The class comprises SDEs that arise from nonlinear transformations of slow-fast stochastic systems. This assumption mimics the situation when we observe the system in terms of the variable $x\in\R^D$ which is an unknown function of a~$D^s$-dimensional slow variable $y$ and a~$D^f$-dimensional fast variable $z$. In this case, the coordinates of $x$ contain both slow and fast dynamics, and we seek to find an observable $\smap$, which we term \emph{the slow map}, such that $\smap(x)=\phi(y)$, i.e., the values of $\smap$ reproduce the slow variable $y$.

In Section~\ref{sec:obs_hid}, we define the observed and hidden systems that give the template for all test cases in our numerical examples and in Section~\ref{sec:features_obs} we delineate three features of the observed process that form the basis for training and testing the networks. In Section~\ref{sec:example}, we give a simple two-dimensional example of the observed-hidden system of interest. Finally, in Section~\ref{sec:slow_map}, we return to the notion of the slow variable from Section~\ref{sec:slow_observables} and use ideas introduced there to give a characterization of the slow map that enables us to assess the accuracy of the trained models.

\subsection{Observed and hidden processes}
\label{sec:obs_hid}

\paragraph{Observed process}
We assume that the observed process $X_t$ is a time-homogeneous diffusion in $\R^{\!D}$, driven by $M$-dimensional noise. More precisely, $X_t$ satisfies the SDE
\begin{equation}\label{eq:obs_sde}
    \der{X_t} = \mu(X_t)\der{t} + \nu(X_t)\der{W_t},
\end{equation}
where $\mu(X_t)$ is the $D\times1$ drift vector, $\nu(X_t)$ is a $D\times M$ diffusion matrix, and $W_t$ is a $M\times1$ vector of independent Brownian motions.
The observed process $X_t$ does not have any obvious slow-fast splitting: the slow and fast modes are mixed together.
However, we suppose that there is an unknown, non-linear observation function~$x=f(y,z)$, where $y$ is a $D^s$-dimensional slow variable and $z$ a $D^f$-dimensional fast variable, with $D^s+D^f=D$, and two corresponding processes $Y_t$ and $Z_t$ such that
\begin{equation}
    X_t = f(Y_t, Z_t).
\end{equation}

\paragraph{Hidden process}
Let us consider a small parameter $\ep\ll1$ that represents the time-scale separation. The hidden process $(Y_t,Z_t)$ is governed by the following slow-fast system of SDEs

\begin{gather}\label{eq:sf_sde}
    \begin{aligned}
        \der{Y_t} &= \mu^s(Y_t, Z_t)\der{t} + \nu^s(Y_t, Z_t)\der{U_t},\\[1em]
        \der{Z_t} &= \frac{\mu^f(Y_t, Z_t)}{\ep}\der{t} + \frac{\nu^f(Y_t, Z_t)}{\sqrt{\ep}}\der{V_t}.
    \end{aligned}
\end{gather}

In~\eqref{eq:sf_sde}, the statistics of the fast variable $z$ quickly becomes tied to the value of the slow variable $y$. We formalize this property by assuming that for fixed $y$ the process $Z_t\cbar y$ that solves the equation
\begin{equation*}
    \der{Z_t} = \frac{\mu^f(y, Z_t)}{\ep}\der{t} + \frac{\nu^f(y, Z_t)}{\sqrt{\ep}}\der{V_t}
\end{equation*}
has an invariant measure $\invm_{\!y}$. The paths of $(Y_t,Z_t)$ fluctuate rapidly along hyper-surfaces $y=\mrm{const}$ and the process $Z_t\cbar y$ equilibrates quickly to $\invm_{\!y}$. Under certain stability assumptions and in the small noise case, it can be demonstrated that the trajectories of~\eqref{eq:sf_sde} are concentrated in a ``layer'' surrounding a manifold of the form
\begin{equation*}
    z(y) = \{\mu^f(y,z) = 0\} +\bo(\ep)
\end{equation*}
and there is a slow spread in the $y$-direction~\cite{berglund_geometric_2003}. We assume that such stability holds in our case, i.e.~we do not consider systems with noise induced transitions, and that the value $z(y)$ is given by the mean of the quasi-invariant distribution $\invm_{\!y}$. The ratio between the speed of fluctuations in the $z$-direction and the evolution along the manifold $z(y)$ is given by $1/\ep$. The correct slow map for $(Y_t, Z_t)$ is simply the projection onto $y$ (as illustrated in Figure~\ref{fig:sin2d_pat_hid}). 

\subsection{Multiscale features of the observed system}
\label{sec:features_obs}
The multiscale nature of the observed process $X_t$ comes from it being a nonlinear transformation of a hidden slow-fast system~\eqref{eq:sf_sde}. Though, in practice, we do not have access to the underlying slow-fast process, there are three main features that arise in the dynamics of $X_t$ due to the presence of the hidden system. These features are ``computationally accessible'', as long as we know the observed SDE~\eqref{eq:obs_sde} explicitly or are able to freely simulate the trajectories of $X_t$.

\begin{enumerate}
    \item \textbf{Local slow and fast directions}
    
    These directions are encoded in the spectral decomposition of the (variance-)covariance matrix of the observed process
    \begin{equation}\label{eq:var_cov}
        \si(x) = \nu(x)\tp{\nu(x)}.
    \end{equation}
     For two-time-scale slow-fast hidden systems the spectrum of $\si(x)$ is divided for each $x$ into two clusters of eigenvalues separated by a large gap.
    The eigenvectors corresponding to the largest eigenvalues identify directions of large noise variability whereas the ones that correspond to small eigenvalues span the local slow subspace.
    
    \item \textbf{Fast fibers}
    
    The fast fibers form a family $\{\ffib_y:\ y\in\R^{D^s}\}$ of $D^f$-dimensional submanifolds of $\R^D$ along which the fast dynamics occurs. The fast fibers foliate the $D$-dimensional state space of the observed system and are tangential to the local fast directions. This family arises from the corresponding $D^f$-dimensional fast hyper-planes of the hidden slow-fast process~\eqref{eq:sf_sde}, which are  aligned with the direction of the fast variable $z$.
    Moreover, by pushing forward the quasi-invariant measures from the fast hyper-planes of~\eqref{eq:sf_sde}, we obtain the corresponding family of probability measures $\invm_{\!y}$ on the fibers $\ffib_y$.
    
    \item \textbf{Slow manifold}
    
    The slow manifold $\sman$ is a $D^s$-dimensional submanifold of $\R^D$ transversal to the fast fibers and such that the slow evolution of process $X_t$ progresses along $\sman$. It corresponds to the adiabatic manifold $z(y)$ of the hidden system discussed in Section~\ref{sec:obs_hid}. The process~$X_t$ does not necessarily evolve on or close to the manifold $\sman$, due to random fluctuations in the transverse fast directions, but for any $x\in\sman$, the statistics of $X^x_t = \Exp[][X_t\cbar X_0=x]$ for a short time windows are fully determined by $x$.
\end{enumerate}
In the subsequent sections, we discuss some numerical approaches for approximating these features. In particular, owing to the stability assumptions on the hidden slow-fast system from Section~\ref{sec:obs_hid}, the points on the slow manifold $\sman$ can be approximated by averaging a sample from the quasi-invariant measures $\invm_{\!y}$. This property provides the basis for computations of the values of the projection onto $\sman$ presented in Section~\ref{sec:dataset}.

\subsection{Example: two dimensional SDE with periodic slow variable}
\label{sec:example}
As the first illustrative example, let us examine a set of two coupled equations from~\cite{froyland_trajectory-free_2016}. The hidden slow-fast system, dating back to~\cite{crommelin_reconstruction_2006}, reads
\begin{gather}\label{eq:sin_hid}
    \begin{aligned}
        \der{Y_t} &= \sin(Z_t)\der{t} + \sqrt{1+\frac{1}{2}\sin(Z_t)}\der{U_t},\\[1em]
        \der{Z_t} &= \frac{\sin(Y_t) - Z_t}{\ep}\der{t} + \frac{1}{\sqrt{\ep}}\der{V_t},
    \end{aligned}
\end{gather}
and is defined for $(y,z)\in[0,2\pi]\times\R$ with periodic boundary conditions on $y$. When we fix $y$, the fast variable $z$ follows the Ornstein-Uhlenbeck process with mean $\sin(y)$ and variance $1/2$. The slow variable $y$ itself moves close to the manifold given by equation $z=\sin(y)$. For all simulations, we set $\ep=10^{-3}$.

\begin{figure}
    \centering
    \includegraphics[draft=false,width=\textwidth]{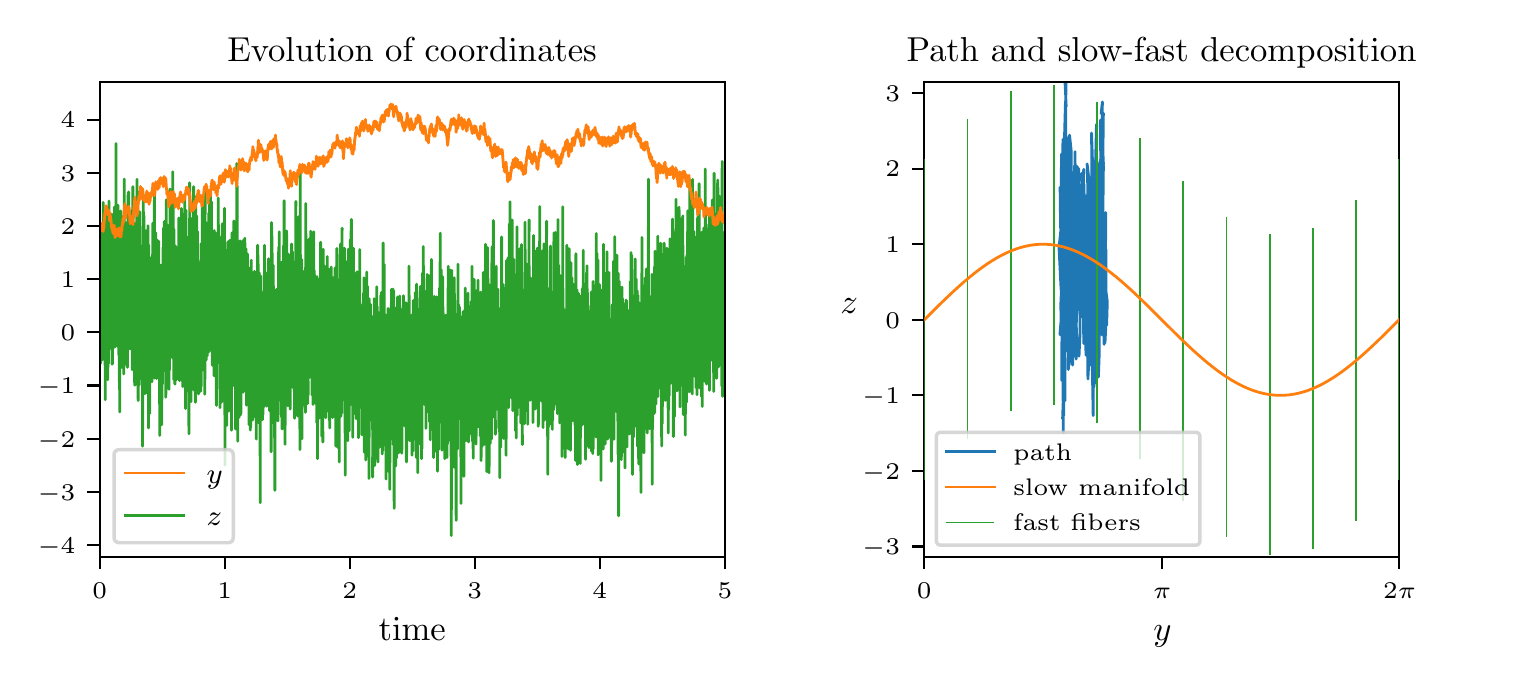}
    \caption{Visualization of a single path of~\eqref{eq:sin_hid}. (Left) The evolution of individual coordinates of the path. The difference between variability of $z$ and $y$ on short time intervals illustrates the time-scale separation between the two variables. (Right) The initial portion of the simulated path as seen in the $(y,z)$ state space (blue) with the slow manifold of the system (orange curve) and a few sections of fast fibers (green vertical segments) superimposed. The path fluctuates rapidly along the fast fibers while slowly drifting to the right such that its mean stays on the slow manifold.}
    \label{fig:sin2d_pat_hid}
\end{figure}

In Figure~\ref{fig:sin2d_pat_hid}, we display a sample path of the slow-fast system~\eqref{eq:sin_hid} with the individual coordinates in the left panel and the systems trajectory in the state space, superimposed by fast fibers and slow manifold, on the right. The slow manifold is given by the relation $z = \sin(x)$ that we obtain by rewriting the fast equation in~\eqref{eq:sin_hid} as $\ep\der{Z}_t = [\sin(Y_t)-Z_t]\der{t} +\sqrt{\ep}\der{V}_t$ and setting $\ep=0$. The fast fibers align with the $z$ axis due to the separation of slow and fast modes in~\eqref{eq:sin_hid}. This separation is also visible in the difference between the variability of the coordinates: the fast variable $z$ fluctuates much more rapidly then the slow $y$ (see the left panel).

\begin{figure}
    \centering
    \includegraphics[draft=false,width=\textwidth]{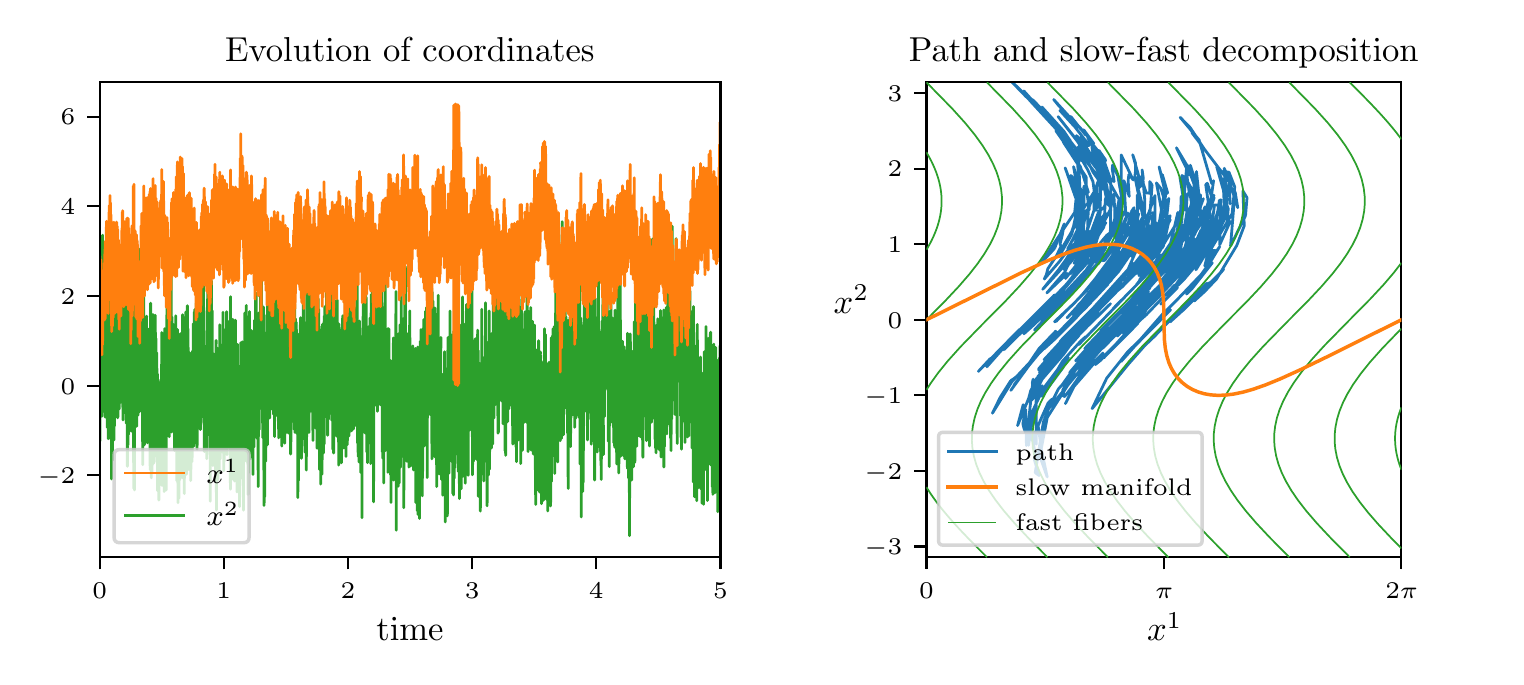}
    \caption{Visualization of a single path of~\eqref{eq:sin_obs}. (Left) The evolution of the individual coordinates from the path. In contrast to Figure~\ref{fig:sin2d_pat_hid}, there is no clear difference between the variability of the coordinates on any time-scale. (Right) The initial portion of the simulated path (blue) as seen in the $(x^1,x^2)$ state space with the slow manifold of the system (orange curve) and a few sections of fast fibers (green sinusoidal curves) superimposed. We observe the same behavior as in Figure~\ref{fig:sin2d_pat_hid} although the geometry of the fast fibers is more complex due to non-linear mixing of slow and fast variables.}
    \label{fig:sin2d_pat_obs}
\end{figure}

To obtain the observed system, we consider the hidden system after a change of coordinates
\begin{equation}\label{eq:sin_transf}
f\colon(y, z)\mapsto (x^1,x^2)=(y+\sin(z), z).
\end{equation}
The mapping $f$ is indeed a diffeomorphism with an inverse that takes $(x^1,x^2)$ to $(x^1-\sin(x^2), x^2)$. After applying~\eqref{eq:sin_transf} to~\eqref{eq:sin_hid}, the resulting observed process $X_t=(X^1_t,X^2_t)$ follows the SDE
\begin{gather}\label{eq:sin_obs}
    \begin{aligned}
        \der{X^1_t} &= \bigg[\sin(X^2_t)+\frac{\cos(X^2_t)}{\ep}\big(\sin(X^1_t-\sin(X^2_t)\big)-\frac{\sin(X^2_t)}{2\ep}\bigg]\der{t}\\[.5em]
        &\hspace{13.5em}+ \sqrt{1+\frac{\sin(X^2_t)}{2}}\der{U_t} + \frac{1}{\sqrt{\ep}}\cos(X^2_t)\der{V}_t,\\[.5em]
        \der{X^2_t} &= \frac{\sin\big(X^1_t-\sin(X^2_t)\big) - X^2_t}{\ep}\der{t} + \frac{1}{\sqrt{\ep}}\der{V_t}.
    \end{aligned}
\end{gather}
Note that now the terms involving $\ep$ appear in the formulas for both coordinates. In the left panel of Figure~\ref{fig:sin2d_pat_obs}, we demonstrate that indeed both coordinates of process~\eqref{eq:sin_obs} fluctuate on the same time-scale; there is no separation present in the coordinates of the observed system. Having the change of coordinates~\eqref{eq:sin_transf} explicitly available, we can recover the observed slow manifold and the fast fibers as images through this transformation of respective submanifolds of the state space of the hidden system. In the right panel of~Figure~\ref{fig:sin2d_pat_obs}, we observe essentially the same behavior as in the case of slow-fast SDE~\eqref{eq:sin_hid}, albeit more complicated: the process fluctuates rapidly along fast fibers while its mean drifts on the slow manifold. That the fast fibers are no longer aligned with any axis reflects the mixing of slow and fast modes in~\eqref{eq:sin_obs}.

\begin{figure}
    \centering
    \includegraphics[draft=false,width=\textwidth]{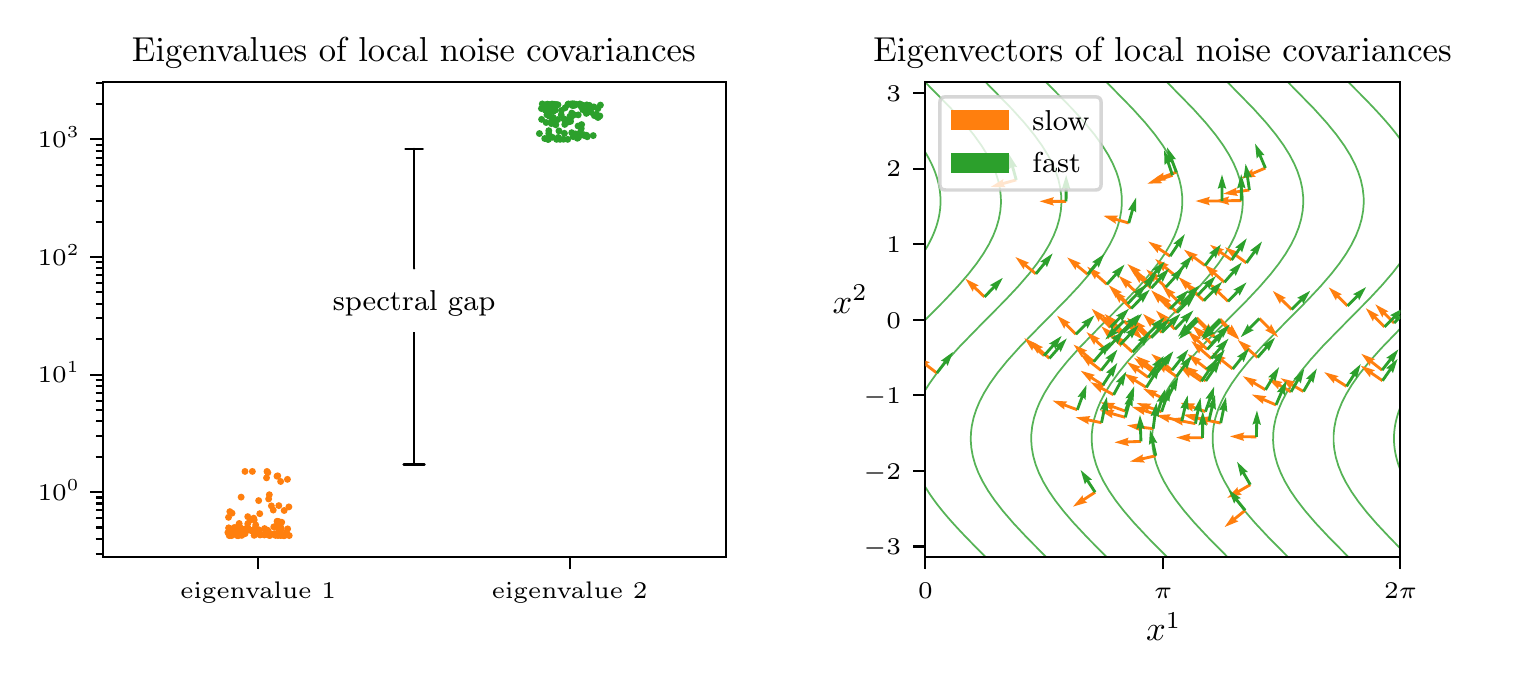}
    \caption{Eigen-decompositions of local noise covariances $\si(x)$ of the observed process~\eqref{eq:sin_obs} evaluated on 100 data points from a sample path. (Left) The categorical scatter plot of the eigenvalues of $\si(x)$ based on their magnitude. The eigenvalues cluster into two groups and the width of the spectral gap corresponds roughly to the time-scale separation equal $1/\ep = 10^3$. (Right) The eigenvectors of $\si(x)$ colored correspondingly to their associated eigenvalues. The eigenvectors corresponding to the large (fast) eigenvalues are tangential to the fast fibers (green curves).}
    \label{fig:sin2d_lnc_eigvals}
\end{figure}

Finally, let us look into the information contained in the covariance matrices $\si(x)$ of the observed process~\eqref{eq:sin_obs}. For this we select 100 points $x_n$ from the simulated path of the observed process and perform an eigen-decomposition of the associated covariance matrices $\si(x_n)$. We present the results in Figure~\ref{fig:sin2d_lnc_eigvals}.

We can detect the time scale separation by looking at the eigenvalues of local noise covariance matrices at each data point. In the left panel of Figure~\ref{fig:sin2d_lnc_eigvals}, we plot, for each data point $x_n$, the two eigenvalues of the covariance matrix $\si(x_n)$. The spectral gap of size roughly $10^3$, equal to the original time-scale separation $1/\ep$, is clearly visible. In the right panel, we visualize the slow and fast local directions at $x_n$ given by the normalized eigenvectors of $\si(x_n)$. Observe that the fast directions, instead of being aligned with one axis, now follow a spatial sinusoidal pattern resulting from the non-linear transformation~\eqref{eq:sin_transf} and they are tangential to the fast fibers of the observed system. Therefore, even when we do not know the underlying transformation from the hidden to the observed system, we can recover the underlying structure of fibers from the eigenvalues of local noise covariances. We shall use this insight to assess the trained models in Section~\ref{sec:ortho_err}.

\subsection{Characterizing the slow map}\label{sec:slow_map}
We combine certain features of the observed process with our tentative definition of slow observable to obtain a characterization of the slow map. This characterization does not require the knowledge of the hidden system nor the underlying transformation from the hidden to the observed process, thus making it a viable approach on which we will build a general computational procedure.

As discussed in Section~\ref{sec:features_obs}, for any observed system~\eqref{eq:obs_sde}, we have a foliation $\{\ffib_y:\ y\in\R^{D^s}\}$ of $\R^D$ into $D^f$-dimensional fast fibers $\ffib_y$.
Intuitively, we expect that a slow map, by exploring the correlations of fast local processes, will collapse the fast fibers of the observed process. Following this intuition, we posit that a slow map is any mapping $\smap\from\R^D\to\R^{D^s}$ that satisfies
\begin{center}
    for every $y\in\img\smap$ there exists a unique fast fiber $\ffib_{y'}$ such that $\levs_y \approx \ffib_{y'}$,
\end{center}
where $\levs_y = \{x\in\R^{D}:\ \smap(x)=y\}$ is the $y$-level set of $\smap$. Therefore, we expect the level sets of a slow map to approximate  the fast fibers of the observed system. If this is the case, the local fast directions of the observed process become tangential to the level sets of $\smap$, thus they get compressed the most under the action of~$\smap$. This indicates that $\smap$ will satisfy condition~\eqref{eq:slowobs_condition} whenever $X_t$ is the observed process~\eqref{eq:obs_sde}.

\begin{rem}[Non-uniqueness]\label{rem:non-uniq}
Note that the condition given above does not guarantee uniqueness of $\smap$. Indeed, whenever $\smap$ is a slow map and $\ph\from\R^{D^s}\to\R^{D^s}$ is one-to-one, the function $\ph\circ\smap$ will have level sets that also align with the fast fibers of the observed system.
\end{rem}

The overarching goal of our method is to encode in the artificial neural network a map~$\smap$ such that the two foliations of the observed state space $\R^D$---the foliation into the fast fibers and the foliation into the level sets of $\smap$---align. Because the solution to this problem cannot be unique, the exact values of $\smap$ do not matter. Therefore, to assess the accuracy of the models, we need to look at the quality of approximation of the fast fibers by $\smap$, see Section~\ref{sec:test_error}.

\begin{rem}[Stochastic and time-dependent slow maps]
In this manuscript, we work with the class of systems that arise as a deterministic transform of slow-fast SDE. This implies that, to extract the slow variables, it suffices to search for a deterministic slow maps. In general, the relation between the observed and the underlying slow-fast system may be more complicated and require a time-dependent and stochastic transform to capture the slow variables, see~\cite{roberts_normal_2008}. Memory terms and randomness can be included in neural networks using special architectures, and we leave it for the future work to explore such approaches.
\end{rem}

\section{A method for learning a slow map}\label{sec:method}
In this section, we present an approach to approximate the slow map of an observed system~\eqref{eq:obs_sde} by an artificial neural network. The method relies on the insight we discussed in the previous section, i.e., to train the network so that the level sets of the encoder mapping overlap with the fast fibers of the observed process.

One can envisage two ways of computationally accessing the fast fibers of the observed process: directly, by sampling them using the ability to simulate the observed process, or indirectly, by looking at the eigenvectors associated with the fast eigenvalues of the covariance matrices (these eigenvectors, as we mentioned, span the tangent spaces of the fibers). Sampling the fast fibers is a challenging task; these are subsets of the state space with measure zero hence the drift and dispersion in the stochastic process force it to quickly escape the fiber containing the starting point. Computing the eigenvectors of the fast eigenvalues is a viable procedure, but it is not clear how obtaining the fibers as integral curves would generalize for higher-dimensional system.

We pursue a different approach that relies on approximating the projection map $\proj$ onto the slow manifold along the fast fibers of the observed process. We use averaging of short bursts of the simulation to produce a reliable approximation of the value of this projection at a given point in state space. By imposing a bottleneck in the network, with size equal to the number of slow variables, we force the network to learn the slow map in this layer. Though the approximation of the slow map clearly falls within the supervised learning techniques, we do not really supervise the representation in the bottleneck. Nevertheless, with a number of test cases presented in Sections~\ref{sec:visu_2D} and~\ref{sec:error_prunning}, we show that we can obtain the correct slow maps in the bottleneck layer.

\subsection{Network architecture}\label{sec:architecture}
We employ a fully-connected, feed-forward network $\net$ with input and output dimensions both equal to~$D$: the~dimension of the state space of process~$X_t$. Additionally, the network contains a bottleneck layer of dimension $D^s$, equal to the number of slow variables in system~\eqref{eq:obs_sde}. In the multiscale systems of interest, $D^s$ will be significantly smaller than $D$ and it can be estimated by analyzing the spectrum of the covariance matrices $\si(x)$.

The architecture of the network is, therefore, the same as for an autoencoder. The main difference between our approach and the standard situation, is that we do not train the network to approximate the identity mapping on $\R^D$. Rather we train it in a supervised manner to approximate the slow projection $\proj$. However, owing to this similarity, we will use the same terminology as employed when describing autoencoders: we divide the network into an encoder $\enc$ and a decoder $\dec$ that are sub-networks before and after the bottleneck, respectively, so that $\net=\dec\circ\enc$. Owing to the interpretation of the bottleneck, we call the bottleneck layer the slow view. Though we do not employ autoencoders, we provide in Section~\ref{sec:autoencoder} an illustration of applying unsupervised approach to our problem to show what can go wrong when using an autoencoder directly.

\begin{figure}
    \centering
    \definecolor{mplblue}{RGB}{31, 119, 180}
\definecolor{mplorange}{RGB}{255, 127, 14}
\definecolor{mplgreen}{RGB}{44, 160, 44}

\begin{tikzpicture}[domain=0.5:5.5, very thick, scale=1.3,
    label/.style={%
        postaction={ decorate,transform shape,
        decoration={ markings, mark=at position .41 with \node #1;}}
    }
    ]

\begin{scope}[thick,smooth,mplgreen,domain=-1:6,
    label/.style={%
        postaction={ decorate,transform shape,
        decoration={ markings, mark=at position .78 with \node #1;}}
    }]
    \clip (0.5,-1.5) rectangle (5.5,2);
    \foreach \y in {-1,...,5}
    {
        \ifnum\y=4
            \draw[label={[below]{\small fast fiber $\ffib$}}] plot ({\y + sin(cos(\x r) r)},\x - \y/2 - 1);
        \else
            \draw plot ({\y + sin(cos(\x r) r)},\x - \y/2 - 1);
        \fi
    }
\end{scope}

\begin{scope}[mplorange]
    \draw[smooth,label={[below]{\small slow manifold $\sman$}}]   plot ({\x, sin(\x r)});
    \draw[->] (8,-1.5) -- (8, 2) node[above] {$\R^{D^s}$};
\end{scope}
\draw (0.5,-1.5) rectangle (5.5,2) node[above left]{observed state space $\R^D$};

\def\yfib{5}
\def\xzer{5}
\def\xone{4.3}
\def\xpro{2.6}
\node (p0) at ({\yfib + sin(cos(\xzer r) r)},\xzer - \yfib/2 - 1) {}; 
\node (p1) at ({\yfib + sin(cos(\xone r) r)},\xone - \yfib/2 - 1) {};
\node (pp) at ({\yfib + sin(cos(\xpro r) r)},\xpro - \yfib/2 - 1) {};
\node (ps) at (8, 0) {};

\def\xco{2.7}
\begin{scope}[mplblue]
  \draw[fill] (p0) circle (1.2pt);
  \draw[fill] (p1) circle (1.2pt);
  \draw[fill] (pp) circle (1.2pt);
  \draw[fill] (ps) circle (1.2pt);
\end{scope}

\draw[->,line cap=round,>=stealth] (p1) .. controls (4, .6) and (3.8, -.2) .. (pp) node[midway,left] {$\proj$};

\begin{scope}[->,>=stealth,thick]
    \draw[shorten >=4pt,dashed] (p0) -- (ps);
    \draw[shorten >=4pt] (p1) -- (ps) node[midway, above] {$\enc$};
    \draw[shorten <=4pt] (ps) -- (pp) node[midway, above] {$\dec$};
\end{scope}


\end{tikzpicture}
    \caption{Comparison of the forward action of the the projection~$\proj$ in the observed system with the forward action of the encoder-decoder network $\net=\dec\circ\enc$. The projection $\proj$ moves points along the fast fibers $\ffib$ of the observed system onto the slow manifold $\sman$. The network $\net$ learns to approximate $\proj$ while passing by the intermediate $D^s$-dimensional representation. Since $\proj$ glues all points on a fast fiber~$\ffib$, the encoder~$\enc$ is forced to associate the same representation to these points. That results in the overlap of the encoder level sets with the fast fibers, a salient property of a slow map enunciated in Section~\ref{sec:slow_map}.}
    \label{fig:enc-dec}
\end{figure}
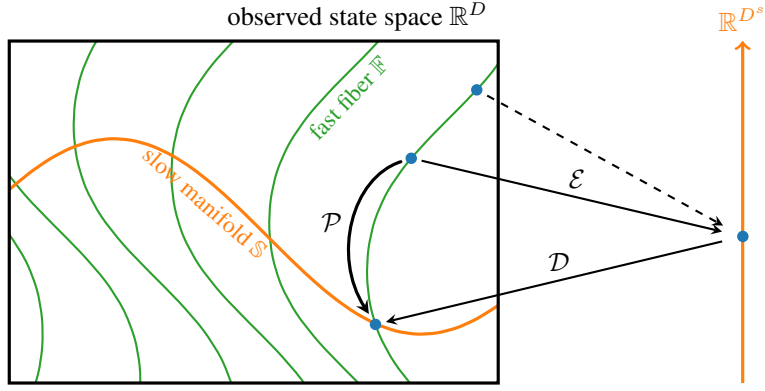

The motivation to use encoder-decoder architecture stems from the nature of the approximated map: the~projection $\proj$. The properties of the slow manifold imply that $\proj$ is an endomorphism of $\R^D$ whose image is a $D^s$-dimensional submanifold of $\R^D$. Consequently, the image of $\proj$ can be parametrized, at least locally, by $D^s$ coordinates. After the training, this parametrization is performed by the decoder $\dec$ while the encoder $\enc$ learns certain $D^s$-dimensional representations of the points of the state space, see Figure~\ref{fig:enc-dec}.

Since $\proj$ collapses the fast fibers, the encoder $\enc$ has to assign all points on the fast fiber to the same representation; with such representation the decoder $\dec$ becomes unaware of the fast fibers. To demonstrate this fact, let us first note that $\enc$ cannot assign the same representation to points on different fast fibers. Otherwise, we could take two points $x,x'\in\R^D$ which have the same representation $\enc(x)=\enc(x')$ yet different projections $\proj(x)\neq\proj(x')$. Then, if the network were to approximate $\proj$, we would have $\net(x)\neq\net(x')$ which contradicts the fact that $\net(x) = \dec(\enc(x)) = \dec(\enc(x'))=\net(x')$. In consequence, the level sets of $\enc$ are $D^f$-dimensional closed submanifolds of the fast fibers $\ffib$ and, since each $\ffib$ is connected, the level sets of the encoder must overlap with the fast fibers.

Therefore, the only way $\net$ can learn the approximation of $\proj$ is when the encoder $\enc$ learns the parametrization of the fast fibers. If this is indeed the case after training the network, for any latent representation $y\in\R^{D^s}$, the $y$-level set of the encoder $\{x\in\R^d:\ \enc(x)=y\}$ will overlap with a fast fiber of the observed system. Hence, $\enc$ will satisfy the main requirement for being a slow map as discussed in Section~\ref{sec:slow_map}.

\subsection{Generating datasets}\label{sec:dataset}
For training and testing the networks we use three ingredients: 1) a sample of data points from the state space that captures the long-time dynamics of the observed process, 2) the corresponding values of the projection onto the slow manifold, and finally 3) the covariance matrices associated with the data points. Therefore, all datasets consist of $M$ instances in the form of triplets
\begin{equation*}
    \big(x_m, \mcl{P}(x_m), \si(x_m)\big),\quad m=1,\dotsc,M,
\end{equation*}
where $x_m$ are the data points in $\R^D$, $\mcl{P}(x_m)$ the corresponding projections of data points onto the slow manifold of observed system~\eqref{eq:obs_sde}, and $\si(x_m)$ the associated local noise covariances. Below, we delineate the exact procedure used to generate our datasets.

\paragraph{Data points}
The data points come from a single trajectory of the observed process $X_t$ sampled at a uniform time interval~$\De t$ over $N$ time steps. To simulate this trajectory, we employ basic Euler-Maruyama method
\begin{equation}\label{eq:em}
    X_{n+1} = X_n + \mu(t_n, X_n)\,\De t + \nu(t_n,X_n)\,\De W_{n+1},\quad X_0 = x,
\end{equation}
where $\De W_{n+1} = W_{t_{n+1}} - W_{t_n}$ are Brownian increments and $x\in\R^D$ is a fixed starting point. For stability reasons,
the time step $\De t$ needs to be a fraction of the fast time scale~$\ep$, hence~$\De t < \ep$. However, to ensure the slow variables had enough time to appreciably evolve, we need to have observations spanning time intervals much larger than $\ep$, thus $N\De t \gg \ep$. This two requirements result in a large $N$, i.e.,~in a large number of points in the sample trajectory. It is generally unnecessary to use all these points to train the networks since we need only a good sample on the time scale of the slow variables. Therefore, to construct the dataset, we subsample the trajectory by uniformly selecting $M$ points from it.

\paragraph{Values of the projection}
We obtain the projection $\mcl{P}(x_m)$ of the data point $x_m$ onto the slow manifold of the observed system by running many parallel short trajectories all starting at $x_m$ and computing the average of their final position. This procedure gives the approximation of the correct values of the projection owing to the connection between the slow manifold and the quasi-invariant measures on the fast fibers $\ffib_y$, discussed in Section~\ref{sec:features_obs}. To attain good approximation, we need to fix an intermediate time-scale $\ta$: sufficiently large to allow the fast variables to sample the invariant distribution on the fast fiber, but short enough to prevent the variables from evolving on the slow manifold (see also Section~\ref{sec:slow_observables}). For observed system~\eqref{eq:obs_sde}, the relevant time scales on which fast processes evolve are approximated by $1/|\la_k|$, where $\la_k$ are the largest eigenvalues of the covariance matrices $\si$. In general, $\ta$ should be fixed as a small multiple of the average of those inverses. With $\ta$ chosen correctly, the end points of the short trajectories sample the invariant distribution of fast variables with slow part of $x_m$ frozen and, due to the stability of the slow manifold, their average approximates the projection $\mcl{P}(x_m)$.

\paragraph{Local noise covariances}
Finally, to compute $\si(x_m)$, we evaluate the dispersion coefficient $\nu(x_m)$ at all data points $x_m$ and obtain the associated covariance via~\eqref{eq:var_cov}. In case~\eqref{eq:obs_sde} is inaccessible directly---because we use an other simulation method, e.g.,~\emph{Gillespie's stochastic simulation algorithm}---we can estimate $\si(x)$, for a fixed $x\in\R^D$, by $J$ parallel one-step trajectories all starting at $x$. The resulting point cloud $\{x_j\}_{1<j<J}$ samples a~Gaussian distribution with covariance matrix $\si(x)$ whose $(d,d')$-entry we estimate by
\begin{equation}
    \si(x)^{d,d'}\approx\frac{1}{\De tJ}\sum_{j=1}^J x_j^dx_j^{d'} - \bar{x}^d\bar{x}^{d'},
\end{equation}
where $\bar{x}$ denotes the average value of the point cloud.

The precise choice of all additional parameters depends on a specific system under study and is a part of the~data preprocessing. Therefore, we do not elaborate more on the procedures required to make the right choice and assume we have data of good quality.

\begin{figure}
    \centering
    \includegraphics[draft=false,width=\textwidth]{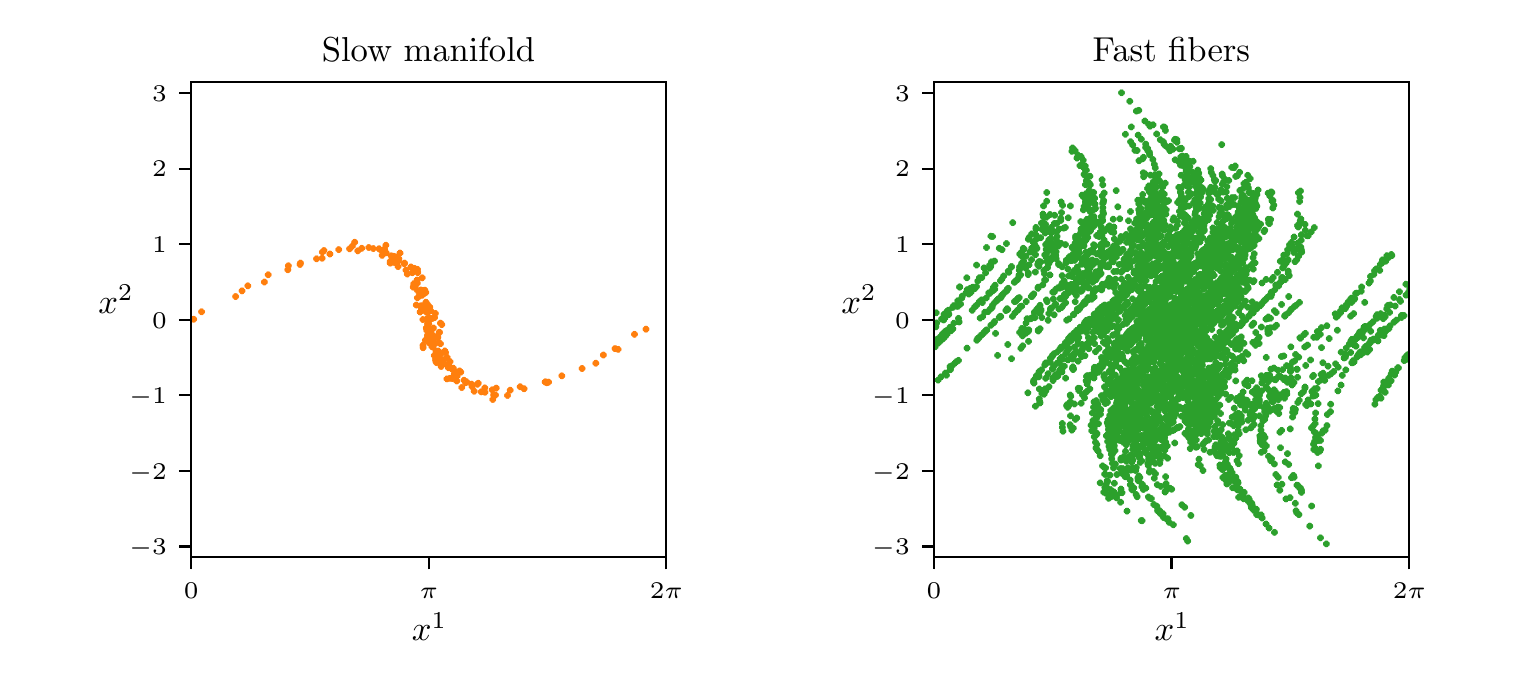}
    \caption{Computing the slow-fast foliation induced by observed process~\eqref{eq:sin_obs}. The slow manifold (left) can be approximated by averaging intermediate bursts of path simulations. See discussion in Section~\ref{sec:dataset} for more details. The fast fibers (right) emerge when we run a very short path simulations starting from a set of data points.}
    \label{fig:sin2d_fibs_sman}
\end{figure}

We present the spectral decomposition of the local noise covariances for system~\eqref{eq:sin_obs} in Figure~\ref{fig:sin2d_lnc_eigvals}. To give an illustration of the two other sets, we look at the results of sampling the slow manifold and the fast fibers for this system.
We plot in Figure~\ref{fig:sin2d_fibs_sman} the approximation of the slow manifold and the fast fibers of~\eqref{eq:sin_obs}. Unlike the corresponding features of the hidden process, these are (computationally) accessible, see Section~\ref{sec:features_obs}, and can be used as additional data attached to the data points sampling the state space.

\subsection{Training}
We train the networks on randomized mini-batches and divide the dataset into training and validation sets. In this phase, as the data instances, we only use the tuples $\big(x_m,\proj(x_m)\big)$, where each $x_m$ serves as the input and~$\proj(x_m)$ as the corresponding label for supervised learning. We always use mean square error loss function.

For each training, we fix a number of maximum epochs and choose the model with the smallest validation loss. We do not employ any stopping criterion, though it is readily available with train-validation splitting. We only monitor the losses during the training to avoid over-fitting and to ensure there is enough epochs for the loss curves to flatten.

\section{Visualization with two-dimensional test systems}\label{sec:visu_2D}
We begin by considering two systems in $\R^2$, each with one slow and one fast variable. This setting allows us to visually compare the slow map with the trained encoder by plotting the values of the slow map against the values of the encoder over a test dataset. Moreover, we can visualize and compare the level sets of both maps to inspect the overlap.

At the end of this section, we present results of training in an unsupervised manner a usual autoencoder to reconstruct the slow manifold. We point out that by doing that, the fast fibers of the observed system are generally not captured correctly.

\subsection{A system with periodic slow variable}\label{sec:sin_analysis}
For the first test, we return to~\eqref{eq:sin_obs}. The slow map, resulting from inverting~\eqref{eq:sin_transf}, reads
\begin{equation}\label{eq:sin_slow}
    \smap:\ x=(x^1, x^2)\mapsto x^1 - \sin(x^2).
\end{equation}

\begin{figure}
    \centering
    \includegraphics[width=\textwidth]{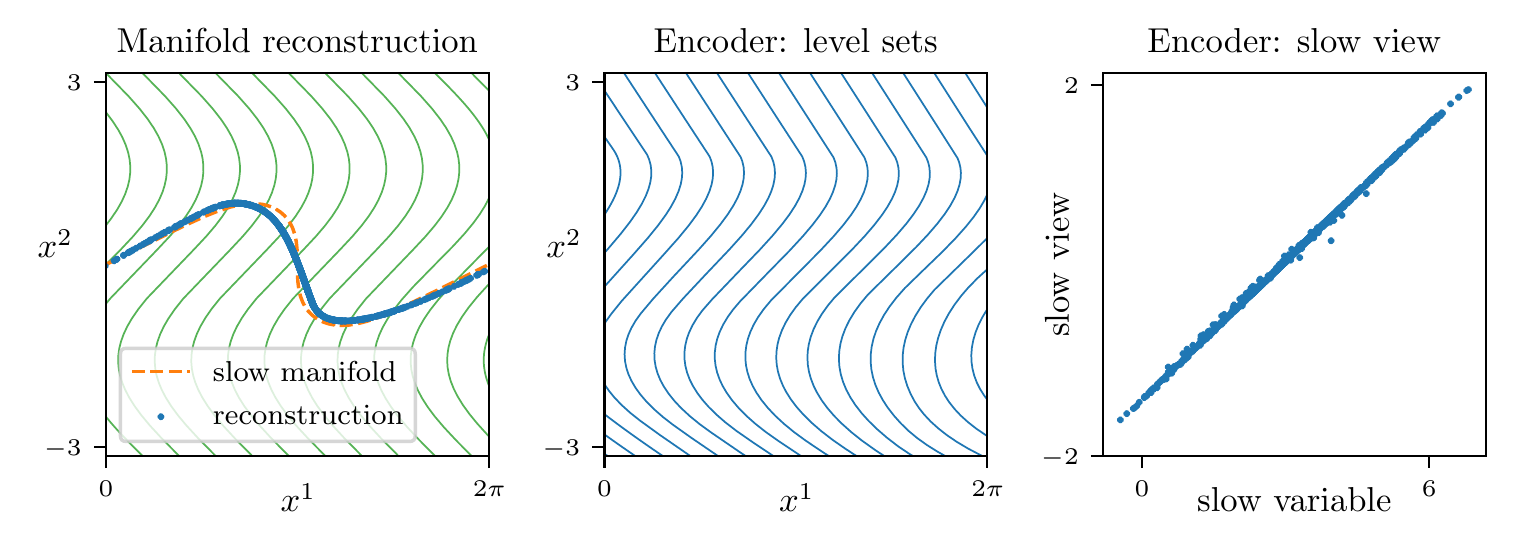}
    \caption{Visualization of the results of trained model for system~\eqref{eq:sin_obs}. (Left) Reconstruction of the values of the projection $\proj$ by the network $\net$ (blue) shows a good accuracy with the underlying slow manifold. (Center) The level sets of the encoder $\enc$. The also show good agreement with the fast fibers of the system displayed in green on the left. (Right) The plot of the values of encoder $\enc$ against the corresponding values of slow variable~\eqref{eq:sin_slow}. The two are clearly related by an affine transformation.}
    \label{fig:sin2d_recon_fibers}
\end{figure}

We trained several models with varying architectures from which we present the smallest one with the best accuracy. We use a training dataset of 1876 instances and validate with 804 instances. The~chosen network has shape 2\,--\,4\,--\,1\,--\,4\,--\,2 with \texttt{ELU} activations and was trained over 3000 epochs with mini batches of size 16 and using \texttt{Adam} optimizer with learning rate equal to 0.003. We plot some crucial features of the trained model and compare them to the known features of~\eqref{eq:sin_obs} in Figure~\ref{fig:sin2d_recon_fibers}. To make a comparison, we compute the values of the network $\net$ and the encoder $\enc$ on the test set of 1321 instances unseen by the network during the training.

In the left panel of Figure~\ref{fig:sin2d_recon_fibers}, we plot in the $(x^1,x^2)$-state space the values of the network $\net$ obtained on the test data (blue). The network is trained to approximate the projection $\proj$ onto the slow manifold $\sman$. Therefore, the network should put any point from the state space onto this manifold, which we refer to as \emph{reconstructing}~$\sman$, and this is confirmed by this plot.

The reconstruction of the slow manifold is only the means to approximating the slow map~$\smap$ by the encoder~$\enc$. In the one-dimensional case, we can directly visualize the accuracy of this approximation by plotting the values of $\enc$ against the corresponding values of the slow map~\eqref{eq:sin_slow}. We display the resulting plot in the right panel. We can see that there is a relation $\enc = \ph\circ\smap$ between the two, where $\ph$ is an affine function.

To gain more insight into the performance of the model, we also display in the central panel of~Figure~\ref{fig:sin2d_recon_fibers} the level sets of the encoder $\enc$. As discussed in Section~\ref{sec:slow_map}, the more the level sets of $\enc$ overlap with the fast fibers of the observed system, equal to the level sets of~\eqref{eq:sin_slow} and plotted in green on the left, the better the approximation. In this case, the model learned the correct sinusoidal shape across the whole range of data.

\subsection{Oscillating half-moons}\label{sec:ohm}
To test the method on a bit more challenging yet still two-dimensional system, we consider the following dynamics from~\cite{singer_detecting_2009}, dubbed the ``oscillating half-moons". The hidden system of SDEs reads
\begin{gather}\label{eq:ohm_hid}
    \begin{aligned}
        \der{Y}_t &= a_1\der{t} + a_2\der{U}_t\\
        \der{Z}_t &= a_3(1-Z_t)\der{t} + a_4\der{V}_t,
    \end{aligned}
\end{gather}
where $a_i$ are constants which we set as $a_1=a_2=10^{-3}$, $a_3=a_4=2.5\cdot 10^{-2}$. The processes $Y_t$ and $Z_t$ are decoupled, with $Y_t$ slowly growing on average with rate $a_1$ and $Z_t$ fluctuating rapidly around~$1$.

\begin{figure}
    \centering
    \includegraphics[width=\textwidth]{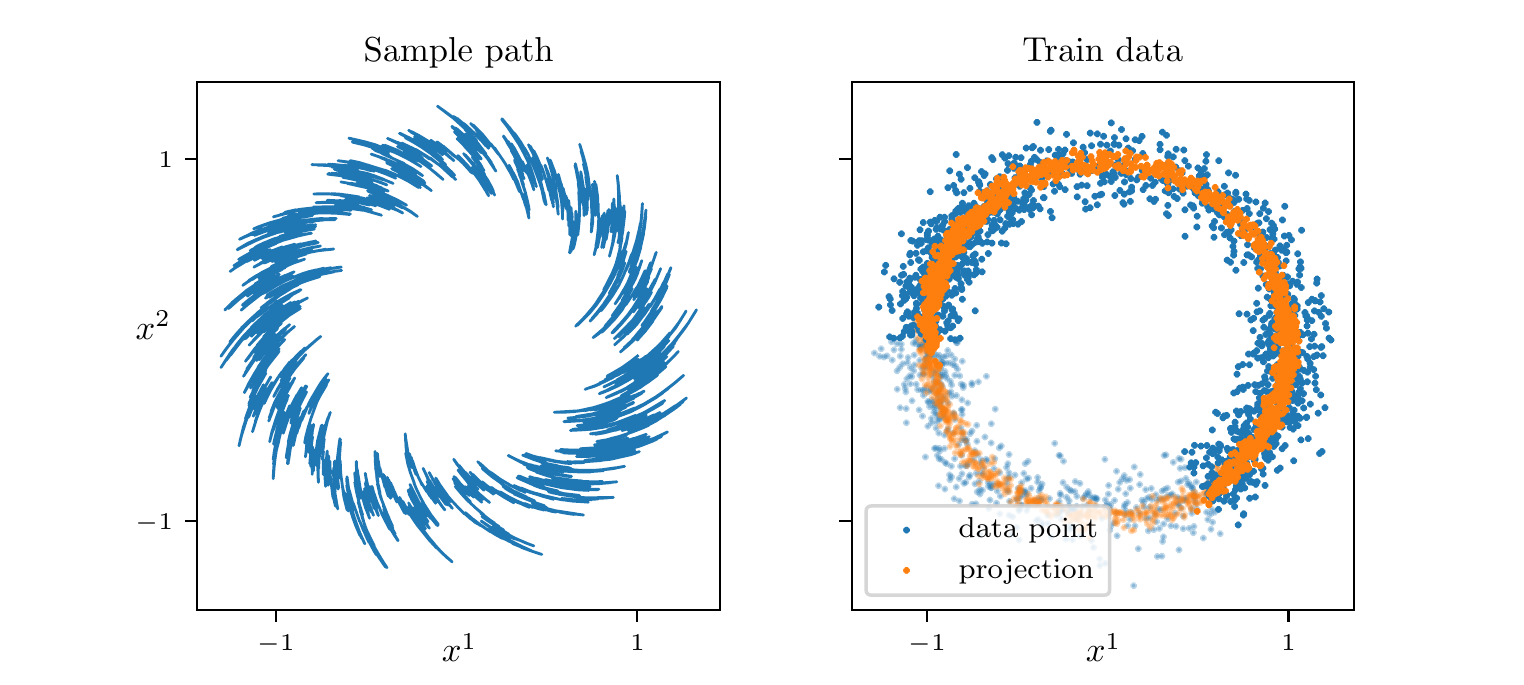}
    \caption{Sample path (left) and training data (right) for the half-moons system from Section~\ref{sec:ohm}. The highlighted data is used first to train the plain encoder-decoder network, for reasons described in the text, whilst we use the full dataset with the modified network from Section~\ref{sec:ohm_polar}.}
    \label{fig:cresc2d_path_data}
\end{figure}

The nonlinear transformation reads
\begin{gather}\label{eq:ohm_transf}
    \begin{aligned}
        x^1 &= z\cos(y+z-1),\\
        x^2 &= z\sin(y+z-1).
    \end{aligned}
\end{gather}
Since $Y_t$ is unstable, this transformation makes the paths of the observed process circle around the center of the observed state space. The fast fluctuations of $Z_t$ are transformed into rapid oscillations along the spirals emanating from the center of the plane. We show an illustrative trajectory in~Figure~\ref{fig:cresc2d_path_data}.

We do not present the full equations of the observed system -- they can be derived (as for all other systems we analyze) from formulas~\eqref{eq:ohm_hid} and~\eqref{eq:ohm_transf} via It\^{o} lemma. Let us only display the slow map for this system:
\begin{equation}\label{eq:ohm_slow}
    y = \arctan(x^2/x^1) + 1 - \sqrt{(x^1)^2 + (x^2)^2},
\end{equation}
which we use henceforth to asses the accuracy of the trained models.

We train the networks on the dataset presented in the right panel of Figure~\ref{fig:cresc2d_path_data}, comprasing of 1600 data points. For the first experiment, we sample only a part of the natural domain of the data points (highlighted points in the right panel, comprising of 1195 data points) in order to avoid the discontinuity in the slow map~\eqref{eq:ohm_slow}. Due to this discontinuity, \eqref{eq:ohm_slow} cannot be accurately represented on the whole $\R^2$ by the continuous encoder map. In Section~\ref{sec:ohm_polar}, we modify the network so that it well approximates the slow map for the complete dataset.

\begin{figure}
    \centering
    \includegraphics[width=\textwidth]{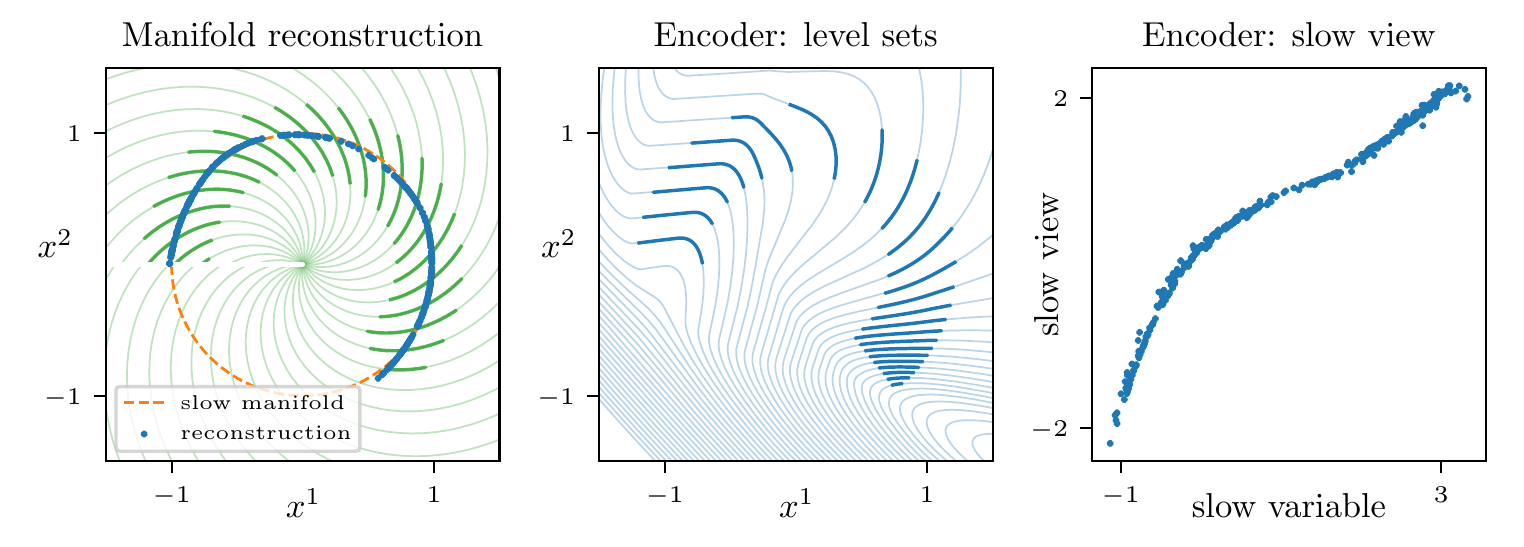}
    \caption{Visualization of the results of trained model for system from Section~\ref{sec:ohm}. (Left) Reconstruction of the slow manifold. (Center) The level sets of the encoder $\enc$. (Right) Values of the encoder $\enc$ plotted agains the corresponding values of the slow map~\eqref{eq:ohm_slow}.}
    \label{fig:cresc2d_recon_fibers}
\end{figure}

First, we train the encoder-decoder of shape 2\,--\,4\,--\,2\,--\,1\,--\,2\,--\,4\,--\,2 and \texttt{ELU} activations over 3000 epochs with mini-batches of size~16 and \texttt{Adam} optimizer with learning rate equal to~0.002. In the leftmost panel of Figure~\ref{fig:cresc2d_recon_fibers}, we can see that the reconstruction of the slow manifold is correct. In the rightmost panel, we can observe that the encoder $\enc$ learned the slow variable accurately. Notice that the relation between the slow map~\eqref{eq:ohm_slow} and $\enc$ is not affine in this case but results from a more complicated transformation.

As noted in Remark~\ref{rem:non-uniq}, the slow map is never unique. To assess the model, we compare the encoder to the slow map that results from inverting the hidden transformation~\eqref{eq:ohm_transf}, which may represent an optimal choice in this way. However, we should not expect that the encoder always finds this particular map, as other choices are readily possible, including affine relation observed in the previous section. Here, we notice only a slight deviation from affinity. This indicates that we found the smallest architecture that can learn the correct slow map, as the more complex networks could stumble upon more complicated, yet correct, encodings.

Looking at the level sets of the encoder in the central panel of Figure~\ref{fig:cresc2d_recon_fibers}, we notice that the global patterns are not the same as for the slow map (left panel, in green). However, if we restrict ourselves to the region where the data is sampled (highlighted lines), the pattern resembles the part of level sets of the slow map in the respective area. As always, we can expect to obtain good accuracy of our models only in the part of the state space where the data resides.


\subsubsection{Incorporating polar coordinates}\label{sec:ohm_polar}
To improve the performance of the network and to achieve better accuracy of the encoder, we can make use of additional knowledge about the system. In the case of oscillating half-moons this is the~periodicity of the~observed dynamics. It is an interesting question whether one can build architectures that discover such properties and automatically adjust. Here, we assume that discovering the periodicity of the system can be achieved without complicated learning approaches and thus used to design the models.  

We observed that the plain encoder-decoder network has difficulties in learning good representation in the slow view due to the discontinuity of the slow variable~\eqref{eq:ohm_slow}. To alleviate this problem, and to demonstrate that one can achieve globally accurate encoding in such a test case, we prepend the encoder $\enc$ with a layer that has no trainable parameters and computes the polar coordinates associated to the inputs it receives. As a result, the remaining part of the encoder, where the parameters are trained, works on transformed data. This approach can be easily generalized to higher-dimensonal systems by making use of (hyper-)spherical coordinates.

We use the network of shape 2\,--\,[2]\,--\,4\,--\,4\,--\,1\,--\,4\,--\,4\,--\,2, where [2] denotes the cartesian-to-polar layer, and the dataset comprising of datapoints sampled from the full path of the system. We train this network over 2000 epochs with mini batches of size 16 and \texttt{Adam} optimizer with learning rate set to~0.002.

\begin{figure}
    \centering
    \includegraphics[width=\textwidth]{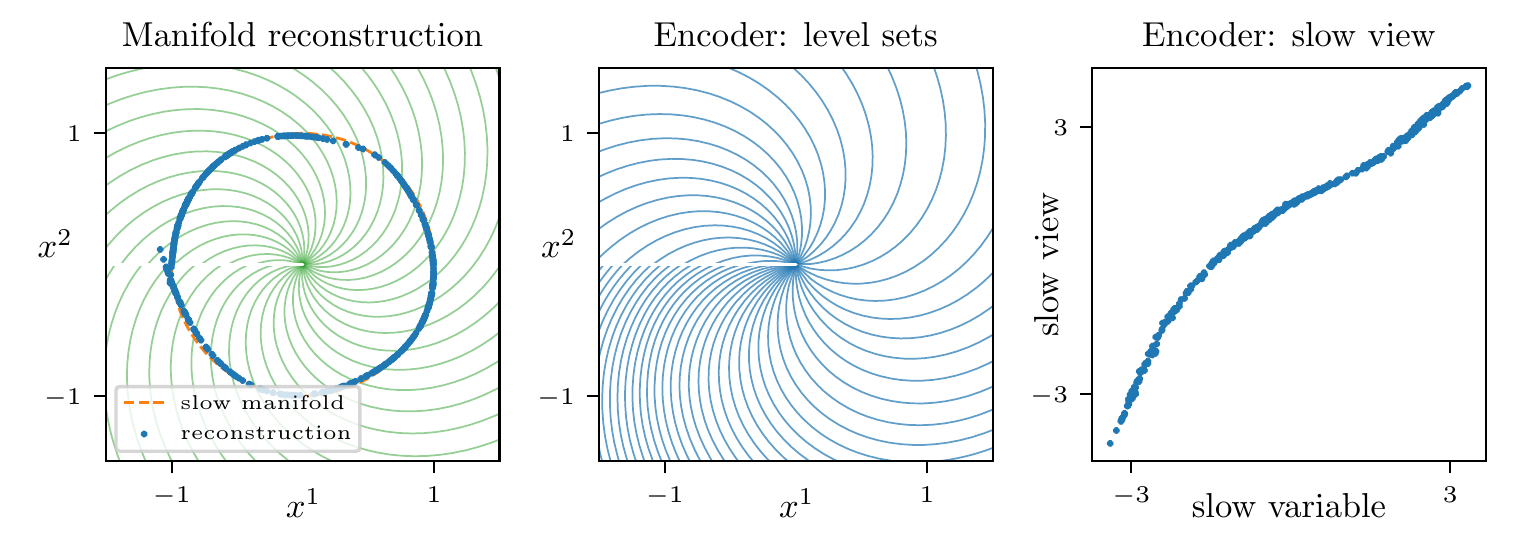}
    \caption{Visualization of the results of trained model for system from Section~\ref{sec:ohm} with encoder prepended with cartesian-to-polar layer. (Left) Reconstruction of the slow manifold. (Center) The level sets of the encoder $\enc$. (Right) Values of the encoder $\enc$ plotted agains the corresponding values of the slow map~\eqref{eq:ohm_slow}.}
    \label{fig:ohm2d_polar_recon_fibers}
\end{figure}

In Figure~\ref{fig:ohm2d_polar_recon_fibers}, we display the recoconstruction of the slow manifold (left), encoder level sets (center) and the accuracy of encoder as compared with~\eqref{eq:ohm_slow} (right). As previously, the reconstruction of the slow manifold is correct but this time the encoder learned the level sets globally and the accuracy of the encoder is good. Notice that the relation between the slow view and the slow variable resembles the one obtained previously, cf.~the right panel of~Figure~\ref{fig:cresc2d_recon_fibers}. This confirms that the previous model, trained on the restricted dataset, agrees on the corresponding region of the state space with the one found in this section.

\subsection{Reconstructing slow manifold via autoencoder}\label{sec:autoencoder}
Before considering higher dimensional systems, let us examine the use of unsupervised learning to parametrize the slow manifold. The goal is to demonstrate that the unsupervised approach based on the dataset of points from a sample path or their projections on the slow manifold do not produce satisfactory results. The method we present in this manuscript combines the two datasets---the points on the sample path as inputs and their respective projections as input labels---to achieve better outcome of learning.

An autoencoder finds a hidden representation of a low-dimensional dataset. The dataset comprising of the points generated by the stochastic system fails, in general, to be a low-dimensional subset of the state space: in the stochastic case the trajectories of the process are only statistically bound to the slow manifolds and can undergo large fluctuations. Therefore, the low-dimensional representation of such a dataset is not well defined and this ambiguity may hinder the discovery of the correct slow parametrization.

In principle, since we can generate the approximation of the slow manifold by performing short bursting coupled with averaging, it is possible to train an autoencoder on the dataset comprising only of  points approximating the slow manifold. This approach allows to build the reduced representation in the latent space, but the question remains as to the relevance of the encoder to the quest of approximating the slow map. The goal of using an autoencoder for parametrizing the manifold is to have a selective reconstruction (only on a subset of the input space) while making the model insensitive to everything outside the manifold. In the cases considered, though, the geometry of the fast fibers is crucial for the problem and may not overlap with the fibers along which the autoencoder projects arbitrary points onto the manifold.

\begin{figure}
    \centering
    \includegraphics[width=\textwidth]{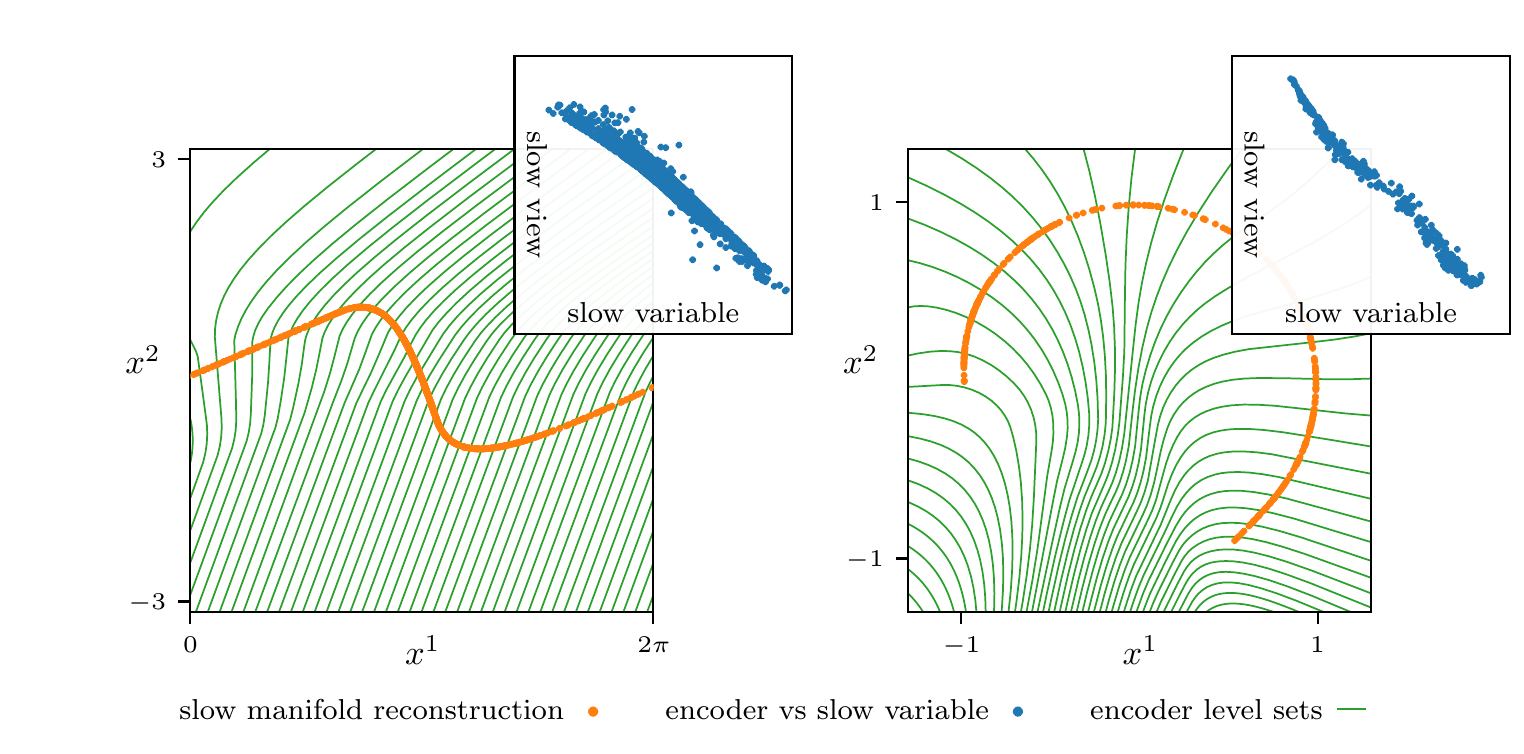}
    \caption{The results of training an autoencoder on the dataset composed solely of the values of the slow projection from Section~\ref{sec:sin_analysis} (left) and~\ref{sec:ohm} (right). Plotted are the reconstruction of the slow manifold by an autoencoder (orange points), the level sets of its encoder (green curves), and the values of encoder agains the slow map (blue points in insets). The reconstruction of the slow manifold is, as expected, accurate, but the level sets of the encoder do not correspond, even locally near the slow manifold, to the fast fibers of the system, cf.~the left panels of Figures~\ref{fig:sin2d_recon_fibers} and~\ref{fig:cresc2d_recon_fibers}. This is reflected in poor correlation between the values of the encoder and the actual slow map as seen by comparing the inset plots with the right panels of Figures~\ref{fig:sin2d_recon_fibers} and~\ref{fig:cresc2d_recon_fibers}.}
    \label{fig:sin2d_cresc2d_auto}
\end{figure}

In Figure~\ref{fig:sin2d_cresc2d_auto}, we demonstrate that we cannot, in general, expect the autoencoder to find the accurate shape of fast fibers, even locally near the slow manifold. For this, we trained two autoencoders on the datasets form Sections~\ref{sec:sin_analysis} and~\ref{sec:ohm}, and display the results in the left and the right panel, respectively. The prediction of the slow manifold is very good (orange). However, the level sets of the encoder (green lines) in both cases do not correspond to the fast fibers of the corresponding system, cf.~Figures~\ref{fig:sin2d_recon_fibers} and~\ref{fig:cresc2d_recon_fibers}. This misalignment, the more pronounced the further from the slow manifold we deviate, indicates that we cannot expect to learn the correct embedding without supervision. This is confirmed in the insets (blue) where we observe that there is poor correlation between the values of the encoder and the values of corresponding slow map on the test sets.

\section{Measuring encoder approximation error and learning essential coordinates}\label{sec:error_prunning}
In this section, we present an approach to measure an error in the accuracy of the encoder as an approximation of a slow map and employ a pruning strategy to select the essential coordinates of the observed process. To illustrate both of these concepts, we apply the proposed method to higher-dimensional test systems, defined in Section~\ref{sec:quad_sde}, for which a simple visualization of the encoder level sets becomes inaccessible.

The error, which we introduce in Section~\ref{sec:test_error}, measures the lack in orthogonality between the gradients of the encoder components and the eigenvectors pointing towards local fast directions. As noted earlier, in order for the encoder to learn an accurate approximation of a slow map, its level sets should globally overlap with the fast fibers of the observed system. Locally, it means that at every data point the row-space of the encoder Jacobian matrix should be orthogonal to the tangent space of the fast fibers at that point. Since we can compute the former using the backpropagation algorithm on the network and the latter is spanned by the eigenvectors accessible through the eigendecomposition of local noise covariance matrix (see Figure~\ref{fig:sin2d_lnc_eigvals}), we obtain a viable error indicator.

We also employ a global pruning of the network as a means to automatically learn the essential coordinates of the observed process. In many high-dimensional systems only a few of the observed variables are involved in the hidden slow dynamics and identifying them constitutes an important task. To address this, we monitor the parameters (weights and biases) of the network during the training and switch off (i.e., fix to zero) the ones with the smallest magnitudes. Such magnitude pruning reduces the complexity of the network and learns the problem-specific architecture along with the parameters. Moreover, and crucial for our application, if the pruning of the parameters of the first layer becomes large enough, it can result in ``cutting off'' of~some input nodes (representing the coordinates of observed process) from the rest of the network. In Section~\ref{sec:essential_coords}, we demonstrate that for a test system in which only a few input nodes build the slow transformation, the remaining nodes are consistently cut off from the network by pruning.

\subsection{Test system with quadratic observation function}\label{sec:quad_sde}
First, let us introduce a model equation we use in this section. It is a slow-fast system that, in addition to the time scale separation $\ep$, is parametrized by the dimensions $D^s$ and $D^f$ of slow and fast dynamics, respectively. This allows to choose different test cases with a given dimensionality and the dimensions of the slow and the fast dynamics set in a particular ratio.

Specifically, we consider a generalization of the system used in~\cite{dsilva_data-driven_2016}. The hidden process is given by a set of uncoupled slow-fast equations
\begin{gather}\label{eq:quad_hid}
    \begin{aligned}
        \der{Y}^d_t &= \der{t} + \der{U}^d_t, & d=1,\dotsc,D^s,\\
        \der{Z}^d_t &= -\frac{1}{\ep}Z^d_t\der{t} + \frac{1}{\sqrt{\ep}}\der{V}^d_t, & d=1,\dotsc,D^f.
    \end{aligned}
\end{gather}
The slow process $Y_t$ is unstable and, as time progresses, it grows on average with constant rate. The fast process $Z_t$ is a stable Ornstein-Uhlenbeck process, so its values center around the origin with normally distributed fluctuations of variance $1/\ep$. In all numerical examples, we set $\ep=10^{-3}$.

We assume that $D^s\leq D^f$ and consider an observed system obtained from~\eqref{eq:quad_hid} by the map
\begin{gather}\label{eq:quad_fun}
    f\ \colon\left\{\ 
    \begin{aligned}
        x^d &= y^d + (z^d)^2, & d=1,\dotsc,D^s,\\
        x^{D^s+d} &= z^d, & d=1,\dotsc,D^f.
    \end{aligned}\right.
\end{gather}
The observation function~\eqref{eq:quad_fun} mixes the squares of the first $D^s$ fast variables with the corresponding slow ones leaving the remaining $D^f-D^s$ fast variables unchanged. These remaining fast variables do not influence the slow dynamics and can be thought of as additional noise present in the system with no impact on the hidden slow process. Clearly, $f$ is a diffeomorphism and the slow map is given by
\begin{equation}\label{eq:quad_slo}
    y^d = x^d - \big(x^{D^s+d}\big)^2,\quad d=1,\dotsc,D^s.
\end{equation}
By the It\^{o} formula, the resulting observed system reads
\begin{gather}\label{eq:quad_obs}
    \begin{aligned}
        \der{X}^d_t &= \frac{1}{\ep}\big[1+\ep- 2\big(X^{D^s+d}_t\big)^2\big]\der{t} &+\ & \der{W}^d_t + \frac{2}{\sqrt{\ep}}X^{D^s+d}_t\der{W}^{D^s+d}_t, & d=1,\dotsc,D^s,\\
        \der{X}^{D^s+d}_t &=-\frac{1}{\ep}X^{D^s+d}_t\der{t} &+\ & \frac{1}{\sqrt{\ep}}\der{W}^{D^s+d}_t, & d=1,\dotsc,D^f.
    \end{aligned}
\end{gather}
As in previous examples, notice that the coefficients involving the time scale separation parameter $\ep$ are mixed in all coordinates of the observed process~\eqref{eq:quad_obs}. However, we can still distinguish two groups of coordinates of size $D^s$ and $D^f$. The first group of $D^s$ coordinates mixes the noise of intensity one with the noise of intensity proportional to $1/\ep$. Additionally, the latter noise correlates the coordinates from the first group with the initial $D^s$ coordinates of the second group, as they are driven by the same Brownian motions. This correlation underpins the hidden slow dynamics in~\eqref{eq:quad_obs} and eventually will be collapsed by a network after the training.

\begin{figure}
    \centering
    \includegraphics[width=\textwidth]{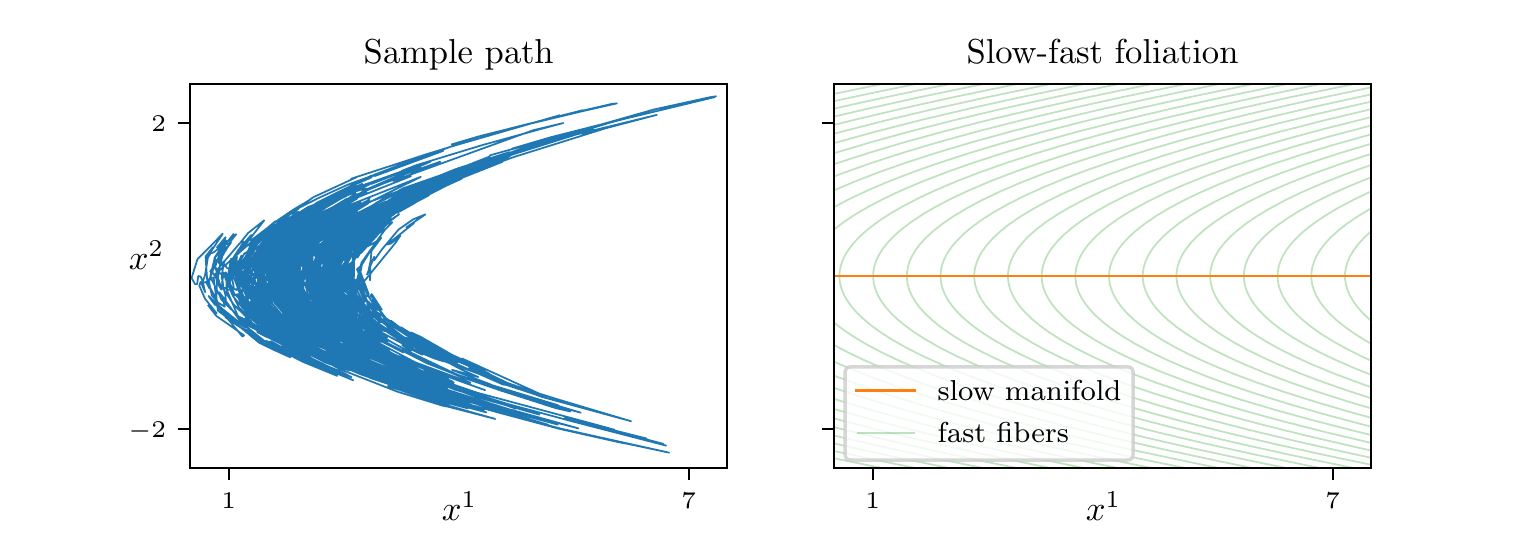}
    \caption{A sample path (left) and the slow manifold with fast fibers (right) of the two-dimensional version of~\eqref{eq:quad_obs}. The parabolic shapes of fast fibers also characterize the higher dimensional versions of~\eqref{eq:quad_obs} when plotting $x^d$ against $x^{D^s+d}$ for any $d=1,\dotsc,D^s$.}
    \label{fig:quad1-1_path}
\end{figure}

A characteristic feature that results from using observation function~\eqref{eq:quad_fun} are the parabolic shapes visible when plotting the $(x^d,x^{D^s+d})$-projections of the paths of the full observed system, illustrated in Figure~\ref{fig:quad1-1_path} for the two-dimensional case. These shapes arise from the intersection of the parabolic fast fibers of the observed process with the $(x^d,x^{D^s+d})$-plane. Finally, let us note that the $D^s$-dimensional slow manifold of~\eqref{eq:quad_obs} is contained in the hyper-plane spanned by the first $D^s$ coordinates.

\subsection{Testing the slow view}\label{sec:test_error}
The goal in this section is to introduce an error measure with which we can assess the accuracy of the slow view without explicitly knowing the slow variable of the observed system. To this end, we consider a four-dimensional version of system~\eqref{eq:quad_obs} with a one-dimensional slow variable ($D^s=1$, $D^f=3$). This setting allows us to visualize, as in the previous section, the accuracy of the encoder by plotting its values on the test dataset against the corresponding values of the slow map~\eqref{eq:quad_slo}. We compare the results of such visualization with the magnitude of the error in Section~\ref{sec:ortho_err} to demonstrate that the error indeed captures the inaccuracies in the slow view. For our experiments, we generate train and test datasets of 2010 and 911 instances, respectively, as described in Section~\ref{sec:dataset}. We use \texttt{ELU} activation function on hidden layers, \texttt{Adam} optimizer with learning rate $0.001$, and train with batches of size $16$. We test three architectures of increasing complexities, see Table~\ref{tab:quad4_architecture}, for which we inspect the reconstruction of the slow manifold and compare the values of the slow view to the slow variable~\eqref{eq:quad_slo}.

\begin{table}
    \centering
    \begin{tabular}{lcrr}
\toprule
 Model   &    Layer sizes    &   Max. epochs &   Min. validation loss \\
\midrule
 Model 1 & 4 - 2 - 1 - 2 - 4 &          2000 &                 0.1063 \\
 Model 2 & 4 - 3 - 1 - 3 - 4 &          2000 &                 0.0029 \\
 Model 3 & 4 - 4 - 1 - 4 - 4 &          2000 &                 0.0029 \\
\bottomrule
\end{tabular}\medskip

    \caption{Architectures and training details for models from Section~\ref{sec:test_error}. Other hyper-parameters remain unchanged across the runs: activation~\texttt{ELU}, \texttt{Adam} optimizer with learning rate~$10^{-3}$, batch size~$16$.}
    \label{tab:quad4_architecture}
\end{table}

\begin{figure}
    \centering
    \includegraphics[width=\textwidth]{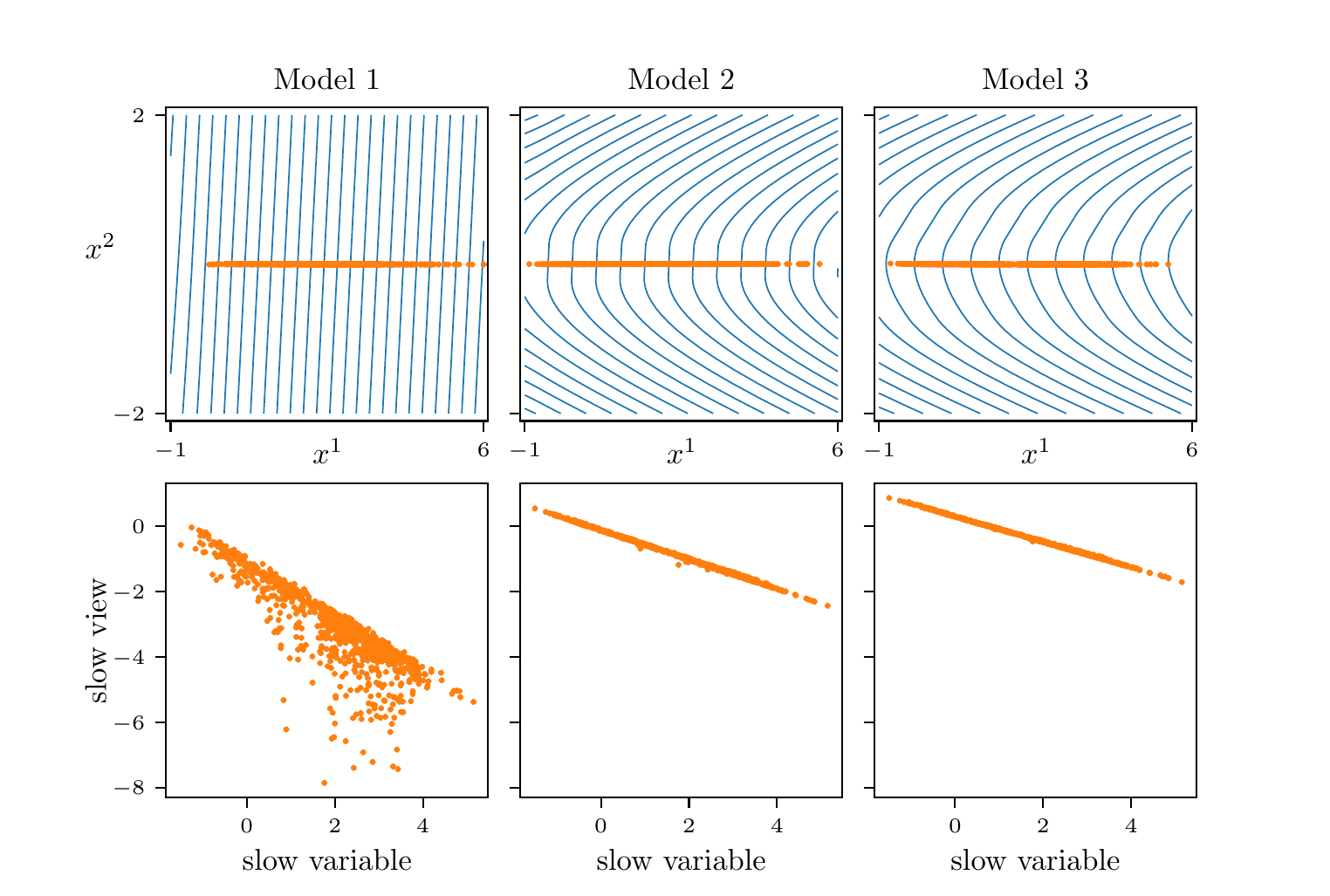}
    \caption{Performance of the models from Section~\ref{sec:test_error}o n a~four-dimensional test dataset for different network architectures (in columns). (Top row) Reconstruction of the slow manifold (orange) and the encoder level sets (blue). (Bottom row) The accuracy of the slow view: values of the encoder versus the corresponding values of the slow map~\eqref{eq:quad_slo}. As the architecture complexity of the networks increases from left to right, the models reach lower validation loss (see Table~\ref{tab:quad4_architecture}) and we get better accuracy of the slow variable (bottom row). Note, however, that the geometry of the slow manifold is well reconstructed by all models (top row, orange), so the difference in accuracy is due mainly to the quality of the approximation of the fast fibers (top row, blue).}
    \label{fig:quad4d_recon_encod_recon}
\end{figure}

Let us first look at the top row of Figure~\ref{fig:quad4d_recon_encod_recon}, where we plot, in the $(x^1,x^2)$-view of the state space, the reconstruction of the slow manifold (orange). The quality of this reconstruction is very similar for three models: all find the correct geometry. However, looking at the value of the smallest validation loss in Table~\ref{tab:quad4_architecture}, we can clearly see that Models~2 and~3 found better minima than Model~1. This discrepancy is due to the simplicity of the geometry of the slow manifold, so that even the smallest architecture, Model~1, can learn the reconstruction correctly, without necessary achieving small loss.

The difference between Model 1 and Models~2 and~3 becomes more pronounced when we look at the approximation of the fast fibers (blue curves) in the top row of Figure~\ref{fig:quad4d_recon_encod_recon}. Only Models~2 and~3 capture the characteristic parabolic shapes of the fast fibers in the $(x^1, x^2)$-plane. Model 1 finds a much simpler encoder level sets, which do not correspond to the fast fibers. As a consequence, the whole network~$\net$ projects data points onto the slow manifold, but not at the correct place on the manifold as determined by the projection~$\proj$. In other words, $\net$ captures the correct geometry of the slow manifold $\sman$ but the approximation of the associated projection mapping~$\proj$ remains insufficient. As a result, the validation loss of Model 1 is much larger than for other two models. This results from the fact that the architecture of Model~1, though sufficient to capture a parametrization of~$\sman$ by the decoder~$\dec$, lacks the~required complexity in the encoder~$\enc$ to reproduce the level sets correctly.

Finally, we can confirm that Models 2 and 3 succeeded in learning the slow variable by looking at the plots in the bottom row of Figure~\ref{fig:quad4d_recon_encod_recon}, which present the scatter plot of slow view of encoder against the slow variable~\eqref{eq:quad_slo}. We observe that, unlike for Model~1, there is an affine relation between the corresponding slow representations of the data points from the test set.

\subsubsection{Measuring local non-orthogonality to the fast directions}\label{sec:ortho_err}
We introduce an error measure that is based on the fact that the slow map's level sets align with the fast fibers of the system. Therefore, at each point in the state space, the subspace spanned by the derivative of the slow map is orthogonal to the tangent space of the fast fibers of the system. Since the encoder should approximate the slow map, one way to measure the quality of this approximation is to measure the lack of orthogonality between the subspace spanned by its derivative and the subspace spanned by the local fast directions.

We can obtain the slow and fast directions by computing the eigenvectors of the covariance matrices, see the right plot in Figure~\ref{fig:sin2d_lnc_eigvals}. The $D^f$ eigenvectors aligning with the local fast directions are the ones corresponding to $D^f$ largest eigenvalues of the local noise covariance matrices. At every data point $x$ for which we want to evaluate the error, we compute the eigendecomposition of $\si(x)$ and assemble the fast eigenvectors into the columns of a $D\times D^f$-dimensional matrix $V^f(x)$.  To obtain the derivative of the encoder with respect to the input $x$ we use backpropagation on $\enc$. This yields the transposed Jacobi matrix $\tp{\jac{\enc}(x)}\in\mat{D}{D^s}$.

Next, we consider the block matrix
\begin{equation*}
    U(x) = \big[V^f(x)\,|\,\tp{\jac{\enc}(x)}\big]\in\mat{D}{D}
\end{equation*}
combining the fast eigenvectors and the derivative of the encoder at the same data point~$x$. The better the encoder approximates the slow map, the more orthogonal the matrix function~$U$ will be across the state space of the system. To measure the lack of orthogonality, we compute the following error
\begin{equation}\label{eq:ortho_err}
    E(x) = \|\tp{U(x)}U(x) - \Id\|_{F} / \sqrt{D},
\end{equation}
where $\Id$ is a $D$-dimensional identity matrix and $\|\,\cdot\,\|_F$ the Frobenius norm.

\begin{figure}
    \centering
    \includegraphics[]{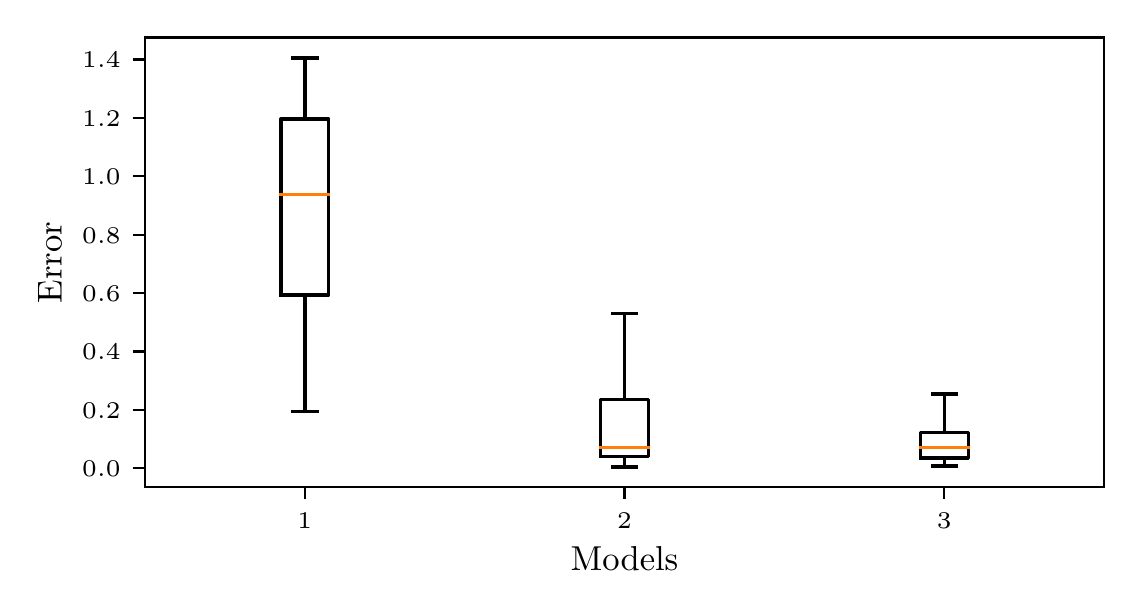}
    \caption{Statistics of the orthogonality error~\eqref{eq:ortho_err} between the encoder gradients and the local fast directions for models from Section~\ref{sec:test_error}:  median  (orange horizontal line), lower and upper quartile (box), 1.5 inter quartile range (whiskers). The statistics is based on the gradients of the encoder computed at all points from the test dataset. In accordance with Figure~\ref{fig:quad4d_recon_encod_recon}, the error decreases as the architecture complexity increases from left to right. We can identify the slightly better accuracy of Model~3 over Model~2, revealed in lower error variability of the former model (see also the discussion in the last paragraph of Section~\ref{sec:ortho_err}.)
    }
    \label{fig:quad4d_enc_err_derivatives}
\end{figure}

In Figure~\ref{fig:quad4d_enc_err_derivatives}, we present the statistics of the error $E(x)$ for Models 1-3, computed over the test dataset. Considering only the median of the error confirms that Model~1 is worse than Models~2 and~3. Moreover, though the medians of Models~2 and~3 are essentially the same, we notice that the variability of the error across the data stays lower for Model~3. This lower variability results from a slightly better representation of the fast fibers close to the slow manifold achieved by Model~3; as seen in Figure~\ref{fig:quad4d_recon_encod_recon}, the level sets for Model~2 have flattened tips of the parabolas. Notice also that the difference between Model~2 and~3 is not reflected in the value of the validation loss (Table~\ref{tab:quad4_architecture}) for these two models.

\subsection{Localizing the essential coordinates via pruning}
\label{sec:essential_coords}
In this section, in addition to learning the slow map by the encoder, we look at the possibility of pinpointing which coordinates of the observed system are involved in the slow variable. It is often the case that in complex, high-dimensional systems only a few of the system's observed coordinates constitute an essential part of the low-dimensional slow dynamics. Finding these coordinates can significantly increase our ability to simulate efficiently the long time dynamics of such system. To accomplish this task, we explore pruning schemes.

Pruning is a method to induce sparsity in the artificial neural network by setting a subset of the model parameters to zero.
While numerous techniques for inducing sparsity have been proposed over the last few years~\cite{gale_state_2019}, in this work we use a simple unstructured iterative pruning that removes small magnitude weights and biases. More specifically, we employ the \href{https://pytorch.org/docs/stable/generated/torch.nn.utils.prune.l1_unstructured.html}{\texttt{l1\_unstructured}} pruning method from the \href{https://pytorch.org/docs/stable/nn.html#utilities}{\texttt{torch.nn.utils}} module. This method removes a fixed amount of (the currently unpruned) parameters with the lowest $L^1$-norm. We aggregate the parameters from a number of specified layers, prior to deciding which ones to prune.

We do not prune the parameters of the bottleneck layer and the output layer. For the remaining layers, the amount of parameters we prune at once varies usually between 3-5\%. Also the pruning schedule---i.e.,~when we start pruning, how often we prune, and when we stop---depends on the model. We usually initialize the pruning strategy after the first 100 epochs and prune every fixed amount of epochs until the desired global sparsity has been reached. After the desired sparsity has been achieved, we continue training for a number of epochs to let the remaining parameters adjust.

We use a ten-dimensional version of~\eqref{eq:quad_obs} with two slow variables as an example test case. According to~\eqref{eq:quad_slo}, the slow map of such system will depend only on the first four coordinates $x^1,\dots,x^4$ of the observed system. The remaining six coordinates do not impact the slow dynamics and constitute noise in the observation.

To test the pruning, we train a number of models with increasingly larger architectures, see Table~\ref{tab:quad10_architecture}, both without and with pruning, on a dataset with $1407$ train instances and $603$ validation instances taken from a sample path as described in Section~\ref{sec:dataset}. We use \texttt{ELU} activation function on hidden layers, \texttt{Adam} optimizer with learning rate $0.001$, and mini-batches of size $16$. The first two models are fairly small and we do not consider the pruned versions. Note also that, besides Model~0, all models reach a similar validation loss. Since we did not observe any over-fitting when training the models, we argue that monitoring solely the validation loss is insufficient to assess the accuracy of the networks trained with our method. This is why other metrics, such as the orthogonality error~\eqref{eq:ortho_err}, should be used in early stopping.

\begin{table}
    \centering
    \begin{tabular}{lcrr}
\toprule
 Model    &              Layer sizes              &   Max. epochs &   Min. validation loss \\
\midrule
 Model 0  &          10 - 5 - 2 - 5 - 10          &          3000 &                 0.0044 \\
 Model 1  &      10 - 6 - 4 - 2 - 4 - 6 - 10      &          3000 &                 0.0029 \\
 Model 2  &      10 - 8 - 4 - 2 - 4 - 8 - 10      &          3500 &                 0.0025 \\
 Model 2p &      10 - 8 - 4 - 2 - 4 - 8 - 10      &          4500 &                 0.0025 \\
 Model 3  &  10 - 8 - 6 - 4 - 2 - 4 - 6 - 8 - 10  &          3500 &                 0.0025 \\
 Model 3p &  10 - 8 - 6 - 4 - 2 - 4 - 6 - 8 - 10  &          5500 &                 0.0024 \\
 Model 4  & 10 - 12 - 8 - 4 - 2 - 4 - 8 - 12 - 10 &          3000 &                 0.0025 \\
 Model 4p & 10 - 12 - 8 - 4 - 2 - 4 - 8 - 12 - 10 &          6000 &                 0.0024 \\
\bottomrule
\end{tabular}\medskip
    \caption{Architectures and training details for models trained on a ten-dimensional dataset where we suffix the pruned models with~`p'. Other hyper-parameters stay the same across the runs: activation function~\texttt{ELU}, \texttt{Adam} optimizer with learning rate~$10^{-3}$, batch size~$16$.}
    \label{tab:quad10_architecture}
\end{table}

\begin{figure}
    \centering
    \includegraphics[width=\textwidth]{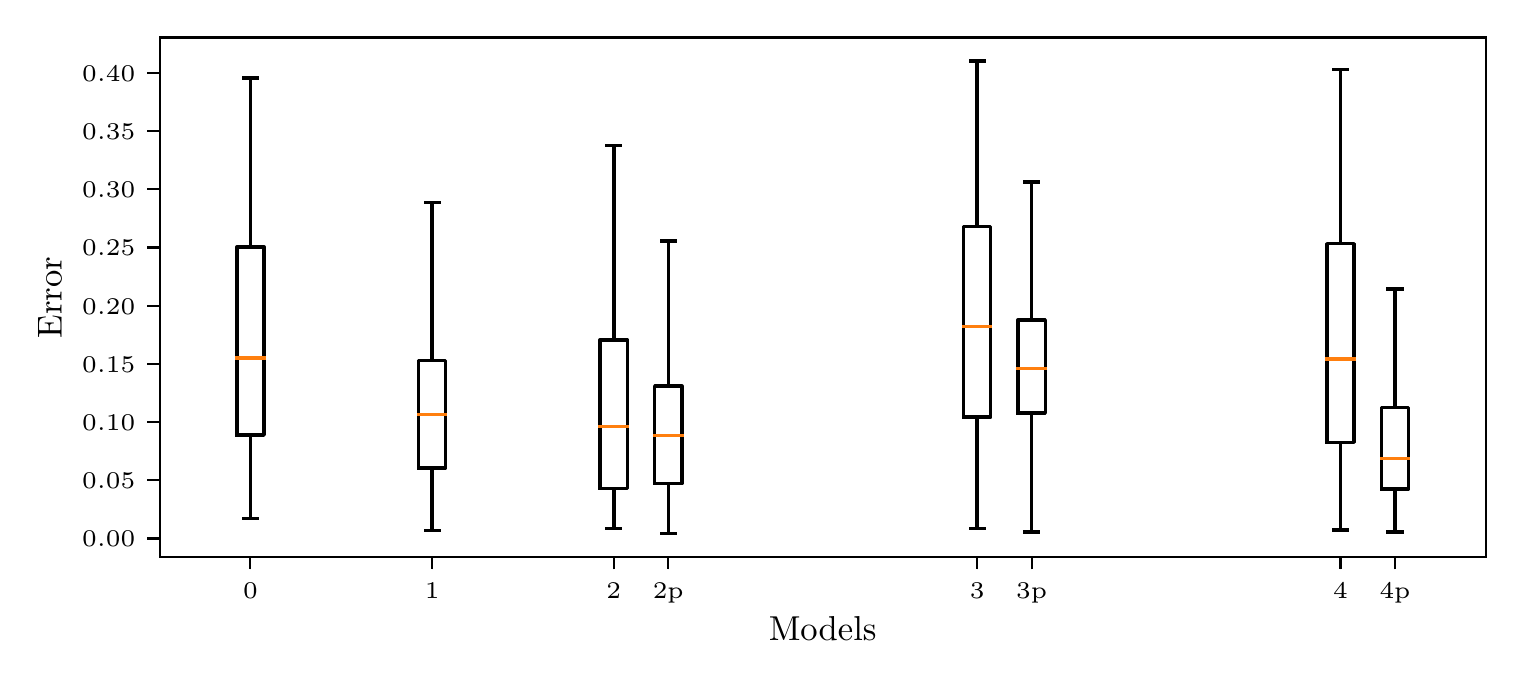}
    \caption{Statistics  of  the  orthogonality of encoder gradients to the local fast directions for  different  models:  median  (horizontal line), lower and upper quartile (box), 1.5 inter quartile range (whiskers).  The plot is based on gradients of the encoder at all points from the test dataset. We suffix the pruned models with~`p'.}
    \label{fig:quad10d_enc_err_derivatives}
\end{figure}

In Figure~\ref{fig:quad10d_enc_err_derivatives} we present the statistics of~\eqref{eq:ortho_err}. We can see that Model~1 performs better than Model~0, in accordance with the decrease of the validation loss. However, as we increase the size of the networks the accuracy of non-pruned models becomes worse: they have larger variability and the average value of the error increases for Models~3 and~4. This is a result of over-fitting: the fibers of the encoder become more complicated than necessary which may worsen the orthogonality. As mentioned in the previous paragraph, this over-fitting is not captured by the validation loss alone.

Adding the pruning strategy to the training process may improve the accuracy of models with architectures that are too rich. Since we do not know beforehand how large the network should be to accurately model the dynamics, removing nodes may prevent over-fitting. This is most visible for the largest model---Model~4. In fact, Model~4p is the most accurate model we could get among the tested architectures. We were also able to diminish the error variability of Model~2 and Model~3.

In general, the proper tuning of the pruning hyper-parameters is challenging and may require substantial testing. We did obtain slight improvements, but this aspect of pruning is not conclusive within the tests we performed. From what we observe by looking at the error, we can also conclude that Model~1 seems to already be a good fit for the test system we use here and not much improvement can be made. To gain a better understanding of the possible improvements in accuracy with pruning strategies a more challenging test case should be used. 

\begin{table}
    \centering
    \begin{tabular}{lcr}
\toprule
 Model    &       Sparsity per layer [\%]        &   Total sparsity [\%] \\
\midrule
 Model 2p &      79 - 19 - 0 - 16 - 10 - 0      &                   30 \\
 Model 3p & 84 - 24 - 32 - 0 - 16 - 43 - 33 - 0 &                   35 \\
 Model 4p & 87 - 65 - 47 - 0 - 25 - 55 - 82 - 0 &                   55 \\
\bottomrule
\end{tabular}\medskip

    \caption{Results of global pruning for each layer of the network. The slow view and output layer were not pruned. Among the remaining layers, the first layer is consistently the most pruned one, indicating that it contains the most redundancy.}
    \label{tab:quad10_sparsity}
\end{table}

\begin{figure}
    \centering
    \includegraphics[width=\textwidth]{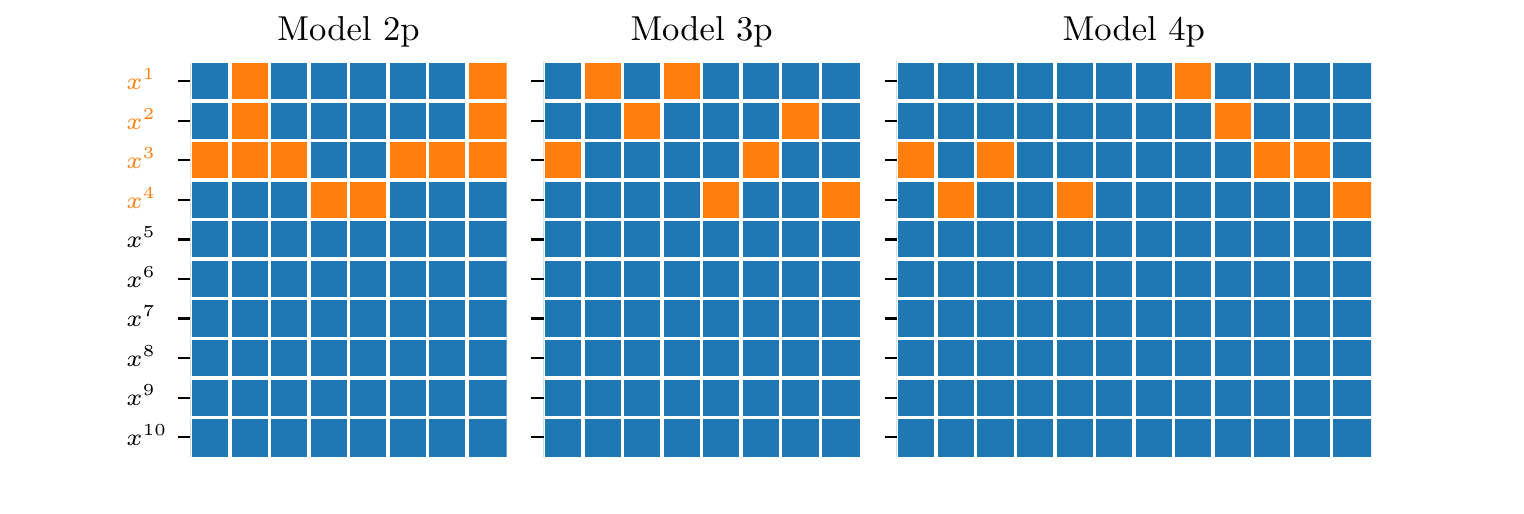}
    \caption{Active (light) and inactive (dark) weights of the first layer of the pruned models. The rows of these matrices represent the weights attached to the coordinates of the observed system. The coordinates with full row inactive are cut off from the rest of the network. For all models, only the first four coordinates, which build up the slow map~\eqref{eq:quad_slo}, remain connected. This indicates that the pruning can discover essential coordinates that constitute the slow variable.}
    \label{fig:quad10d_sparsity}
\end{figure}

We summarize the most important effect of pruning in our case in Table~\ref{tab:quad10_sparsity} and Figure~\ref{fig:quad10d_sparsity}: due to a limited number of observed coordinates involved in the slow dynamics, pruning results in high sparsification of the weights matrix of the first layer of the networks. In Table~\ref{tab:quad10_sparsity}, we compare the amount of sparsification per layer in the pruned networks. Note that we did not turn pruning on for the parameters of the bottleneck layer and the output layer. The biggest percentage of pruned parameters appears consistently in the first layer, although we used a global pruning strategy that looks at the parameters closest to zero across all active layers. This clearly indicates that the first layer is the place where most redundancy is present.

Furthrmore, the pruning cuts off non-essential input nodes from the rest of the network, indicating which input nodes (representing the coordinates of the observed system) influence the latent view. In Figure~\ref{fig:quad10d_sparsity}, we display the activation of weight matrices of the first layers for the pruned models. The orange squares indicate that the weight is different than zero. Notice that rows 5-10 have been completely pruned out: the input nodes corresponding to these rows are no longer connected to the network and thus do not influence the prediction. The rows that still contain active weights correspond to the first four coordinates of the system---exactly those coordinates that constitute the slow dynamics.

\section{Conclusions and outlook}
The main contribution in this paper is to introduce a learning method based on approximation properties of artificial neural networks that is able to discover the slow variables in stochastic systems with multiple time scales. The networks share the architecture with autoencoders, but we train them in a supervised manner to approximate the projection map along the fast fibers onto the slow manifold of the observed system. We demonstrated, on a set of test examples, that the encoder part of the network can successfully learn a correct slow representation. A number of possible improvements in the architecture and training of the networks exists which we did not explore in this initial study but which may produce better results and become useful in applications to more demanding multiscale systems.

We can visually inspect the accuracy of encoders in approximating the slow map on examples with one-dimensional slow representation. To measure this accuracy we introduced an error based on the lack of orthogonality between the gradients of the encoder and the local fast directions. The error can be easily computed using the backpropagation algorithm with respect to the input on the encoder part of the network. This allows to not only compute it for trained models, but also using it as a monitoring metric during training. Since, the accuracy is not always reflected in the value of the loss function, i.e., models with the same validation loss still can have noticeably different accuracies, having a metric based on the error allows to implement early stopping and to improve the results of training.

The architecture of our network is, as for autoencoders, symmetric. However, while the encoder part of the network is responsible for collapsing the fast fibers that foliate the whole state space, the decoder's task is only to parametrize the lower-dimensional slow manifold. Therefore, it seems reasonable to expect the same or better performance using non-symmetric architecture, in which the decoder part is smaller than encoder. On one hand, this can reduce the size of the network, thus speeding up the training. On the other hand, if such architecture indeed fits the problem better, this can prevent the overfitting. Thorough tests of such architectures on various systems could provide evidence for this hypothesis.

Finally, we used pruning to adjust the architecture to a given problem and to prevent overfitting when starting with too large network. Due to the simple nature of our test systems, we only achieved modest improvements in these respects and further tests on more challenging systems should follow. However, with pruning we accomplished a more important goal in our case: we were able to locate the essential coordinates of the observed system that build up the slow map. Knowing such coordinates opens up the possibility to construct more efficient algorithms for the simulation of high-dimensional systems. In this study, we employed the simplest magnitude pruning with aggregation over a number of layers. More intricate techniques are available in the literature and it remains to be seen whether they can perform even better in certain situations. Since our goal is to remove the whole rows from the weight matrix of the first layer, the techniques that impose block structure on the spare weights, more efficient in parallel processing, can provide an interesting modification to explore in the future work.

\section*{Acknowledgments}
This work was partially supported by AFOSR under grant FA9550-17-1-9241.

\printbibliography 

\end{document}